\documentclass[twocolumn]{autart} 
\usepackage{graphicx}      
\usepackage{natbib}        
\usepackage{graphics} 
\usepackage{epsfig} 
\usepackage{amsmath} 
\usepackage{amssymb}  
\usepackage{amsfonts}  
\usepackage{subfig}
\usepackage{bm}
\usepackage{epstopdf}
\usepackage{enumerate}
\usepackage{graphicx} 
\usepackage{picins}
\usepackage{color}
\newcommand{\tcb}{\textcolor{black}}

\newcommand*{\QEDB}{\hfill\ensuremath{\square}}%

\usepackage{blindtext}
\usepackage{url}
\usepackage[hidelinks,bookmarks=false]{hyperref}

\newtheorem{theorem}{Theorem}
\newtheorem{proposition}{Proposition}
\newtheorem{lemma}{Lemma}

\newtheorem{definition}{Definition}
\newtheorem{assumption}{Assumption}
\newtheorem{remark}{Remark}
\newtheorem{example}{Example}

\usepackage{lineno}

\begin{document}
	\begin{frontmatter}
		
		\title{\vspace{-0.3cm}On Composite Foster Functions for a Class of\\ Singularly Perturbed Stochastic Hybrid Inclusions:\\ \tcb{Theory and Applications}\vspace{-0.4cm}} 
		\thanks[footnoteinfo]{Earlier preliminary results of this work were presented at the 2023 IFAC World Congress, Yokohama, Japan.\\
			Research supported in part by AFOSR grant FA9550-22-1-0211 and NSF grant CAREER 2144076.}
		
		\author[a]{Jorge I. Poveda},
		\author[b]{Mahmoud Abdelgalil}
		\address[a]{Department of Electrical and Computer Engineering, University of California, San Diego, La Jolla, 92093, California, USA}
        \address[b]{\tcb{Department of Mechanical and Aerospace Engineering, University at Buffalo, SUNY, Amherst, NY 14260-4400, USA}}
		
		\begin{abstract}               
		We study sufficient conditions for stability and recurrence in a class of singularly perturbed stochastic hybrid dynamical systems. The systems considered combine multi–time-scale deterministic continuous-time dynamics, modeled by constrained differential inclusions, with discrete-time dynamics described by constrained difference inclusions subject to random disturbances. Under suitable regularity assumptions on the dynamics and causality of the associated solutions, we develop a family of composite nonsmooth Lagrange–Foster and Lyapunov–Foster functions that certify stability and recurrence properties by leveraging simpler functions associated with the slow and fast subsystems. Stability is characterized with respect to compact sets, while recurrence is established for bounded open sets. The proposed framework is illustrated through several examples and applications, including the stability analysis of singularly perturbed switching systems with stochastic spontaneous mode transitions, feedback optimization problems with stochastically switching plants, and momentum-based feedback-optimization algorithms with stochastic restarting.
		\end{abstract}
		
		
	\end{frontmatter}
	
	\section{Introduction}

    \vspace{-0.2cm}
    Methods that exploit multiple time scales have long been central to feedback control, with singular perturbation analysis \cite{Krylov1,Kokotovic_SP_Book,Saberi,TeelMoreauNesic2003,Teel96SP} emerging as a foundational tool across applications such as optimal control \cite{Kokotovic_SP_Book}, steady-state optimization \cite{colombino2019online,OnlinePovedaHybrid}, extremum-seeking control \cite{PoTe17Auto}, etc. 
\tcb{More recently, the emergence of hybrid dynamical systems—combining continuous- and discrete-time dynamics—has motivated singular perturbation tools tailored to hybrid settings \cite{Sanfelice:11,abdelgalil2023multi,SP2020Teel,RajebSP,tanwani2024singularly}.} These developments have enabled control and optimization algorithms that explicitly capture interactions between slow and fast hybrid processes, with applications including momentum-based steady-state optimization \cite{OnlinePovedaHybrid}, hybrid extremum-seeking control \cite{PoTe17Auto,abdelgalil2023multi,poveda2025hybrid}, distributed optimization \cite{wang2020stability}, and learning in game-theoretic networks \cite{TAC21Momentum_Nash}. However, many modern systems exhibit inherently stochastic behavior—arising, for example, in stochastic optimization \cite{baradaran2018stochastic}, model-free control under random perturbations \cite{PoTe16CDC}, and synchronization under random graphs \cite{JavedAutomatica20}—which is more naturally modeled within the framework of stochastic hybrid dynamical systems (SHDS) \cite{Hespanha_2007,Teel:14_AutomaticaSurvey}. Despite extensive developments for deterministic ODEs and hybrid systems, singular perturbation theory for stochastic hybrid settings remains relatively underdeveloped. Existing results primarily address specialized cases, such as Markov jump linear systems \cite{HaddadsingularlyPerturbed,Hinfinitesingularly}, finite-horizon formulations \cite{FilarTAC}, or stochastic ODEs, and therefore do not fully capture the general structure of hybrid stochastic dynamics.
%

\vspace{-0.1cm}
This paper introduces a general framework based on composite auxiliary functions for the study of two key properties in singularly perturbed SHDS (SP–SHDS): Uniform Global Asymptotic Stability in Probability (UGASp) and Uniform Global Recurrence (UGR). To characterize these properties in singularly perturbed SHDS, we \tcb{extend} the well-known composite Lyapunov method, studied in \cite{Saberi,tang2025fixed} for ODEs, to construct nonsmooth Lyapunov functions built from simpler functions associated with the system’s slow and fast dynamics. In the context of SHDS, these constructions employ regular Lagrange–Foster and Lyapunov–Foster functions, allowing stability and recurrence to be analyzed using modular components. This approach provides Lyapunov-based certificates for singularly perturbed SHDS with set-valued dynamics. It also unifies and generalizes earlier results obtained for deterministic hybrid systems \cite{Wang2020SP,abdelgalil2023multi} and continuous-time ODEs \cite{Saberi}. The proposed tools leverage the mathematical model of SHDS introduced in \cite{Teel_ANRC}, and are expected to facilitate the design and analysis of feedback, optimization, and estimation algorithms that simultaneously involve stochastic effects and hybrid time-scale behavior—scenarios where robust stability and recurrence guarantees are crucial. 

\vspace{-0.2cm}
To demonstrate the practical relevance of the proposed framework, several examples are presented in the paper, including two representative applications: (a) The stability analysis of singularly perturbed switching systems with stochastic mode transitions, where the deterministic results of \cite{chatterjee2006stability,chatterjee2011stabilizing} are extended by explicitly constructing composite Lyapunov functions from those of the reduced and boundary-layer subsystems; and (b)  real-time steady-state optimization (also called feedback optimization) of dynamical systems using controllers with stochastic momentum restarting, which extends the deterministic smooth and hybrid strategies studied in \cite{colombino2019online,FlorianSteadyState,OnlinePovedaHybrid}, as well as incorporating stochastically switching plants. While stochastic resetting and stochastic restarting techniques in optimization and dynamical systems have been broadly studied in the literature \cite{gupta2022stochastic,belan2018restart,pokutta2020restarting,wang2022scheduled,fu2025hamiltonian}, our theoretical framework demonstrates that analogous stochastic hybrid mechanisms can be employed for the solution of feedback-based optimization problems with plants in the loop. This leads to the development of the first class of stochastic hybrid feedback optimization algorithms that provides probabilistic closed-loop stability guarantees. Finally, we also show how the proposed tools can be used to study feedback optimization problems with stochastically switching plants in the loop, thus extending the results for deterministic switching plants studied in \cite{OnlinePovedaHybrid}.
	
	\vspace{-0.2cm}
	Earlier preliminary results of this work were reported in the conference paper \cite{poveda2023singularly}. The present manuscript significantly extends that contribution by providing the complete proofs, a more comprehensive analysis,  sharper results for recurrence under weaker assumptions, and, additionally, several novel detailed applications in the context of switching systems and feedback optimization. 

 %
	\section{PRELIMINARIES}
	\label{sec:preli}
	\textbf{Notation:} We denote the set of (non negative) real numbers by $(\mathbb{R}_{\geq 0})$ $\mathbb{R}$.
	The set of (non-negative) integers is denoted by $(\mathbb{Z}_{\geq 0})$ $\mathbb{Z}$. Given a 
	set $\mathcal{A} \subset \mathbb{R}^n$ and a vector $z \in \mathbb{R}^n$, we define $|z|_\mathcal{A} := \inf_{y \in \mathcal{A}}\|z-y\|$, and we use $\|\cdot\|$ to denote the standard Euclidean norm. We use $\overline{\mathcal{A}}$ to denote the closure of $\mathcal{A}$, and $
	r\mathbb{B}^\circ$ to denote the open ball (in the Euclidean norm) of appropriate dimension centered around the origin and with radius $r>0$. For ease of notation, given two vectors $u,v \in \mathbb{R}^{n}$, we write $(u,v)$ for $(u^{\top},v^{\top})^{\top}$. Given a matrix $A\in\mathbb{R}^{n_1\times n_2}$, we use $\sigma_{\max}(A)$ to denote the maximum singular value of $A$, and $\lambda_{\max}(A)$ to denote the maximum eigenvalue of $A$. In particular, $\sigma_{\max}(A) = \sqrt{\lambda_{\max}(A A^\top)}$. We also use $\mathrm{Sym}(A)$ to denote the symmetric part of $A$, i.e., $\mathrm{Sym}(A) = \frac{1}{2}(A+A^\top)$. A function $f:\mathbb{R}^n\to\mathbb{R}$ is said to be: a)  $C^k$ if its $k^{th}$ derivative is continuous; and b) radially unbounded if $f(x)\to \infty$ whenever $\|x\|\to\infty$. A function $\alpha:\mathbb{R}_{\geq0}\to\mathbb{R}_{\geq0}$ is said to be: a) of class $\mathcal{G}_{\infty}$ if it is continuous, non-decreasing, and unbounded; b) of class $\mathcal{K}_{\infty}$ if it is zero at zero, continuous, strictly increasing, and unbounded. It is said to be of class $\mathcal{P}_s\mathcal{D}(\mathcal{A})$ if it is continuous and positive semidefinite with respect to $\mathcal{A}$, and of class $\mathcal{P}\mathcal{D}(\mathcal{A})$ when it is continuous and positive definite with respect to $\mathcal{A}$. When $\mathcal{A}=\{0\}$, we simply use $\mathcal{P}_s\mathcal{D}$ and $\mathcal{P}\mathcal{D}$. We use $\overline{\text{co}}(\mathcal{A})$ to denote the closure of the convex hull of the set $\mathcal{A}$, and $\mathbb{I}_{\mathcal{A}}:\mathbb{R}^n\to\{0,1\}$ to denote the standard indicator function. We use $\mathbf{B}(\mathbb{R}^m)$ to denote the Borel $\sigma$-field, and $K\subset\mathbb{R}^m$ is said to be measurable if $K\in\mathbf{B}(\mathbb{R}^m)$.
	
	Let $f:\mathbb{R}^n\to\mathbb{R}$ be a locally Lipschitz function. The Clarke generalized gradient of $f$ at $y\in\text{dom}~f$, is the set $\partial f(y):=\text{co}\{v\in\mathbb{R}^n:\exists y_k\to y,~y_k\notin\mathcal{Z},~\lim_{k\to\infty}\nabla f(y_k)=v\}$, where $\mathcal{Z}$ is the set where the usual gradient $\nabla f$ is not defined, which is of measure zero due to Rademacher's Theorem. We use $\partial_{x_1} f(x_1,x_2)$ and $\partial_{x_2} f(x_1,x_2)$ to denote the partial Clarke gradients. The function $f$ is said to be regular at $y$ if, for every $u\in\mathbb{R}^n$, the directional derivative \tcb{$f'(y;u):=\lim_{s\to 0^+}\frac{f(y+su)-f(y)}{s}$} exists, and the following holds: \tcb{$f'(y;u)=\max\{v^\top u:v\in \partial f(y)\},~\forall~u\in\mathbb{R}^n$}. Typical examples of locally Lipschitz regular functions include $C^1$ and convex functions. A set-valued mapping \tcb{$M:\mathbb{R}^m\rightrightarrows\mathbb{R}^n$} is outer semi-continuous (OSC) if for each $(x_i,y_i)\to(x,y)\in\mathbb{R}^m\times\mathbb{R}^n$ satisfying $y_i\in M(x_i)$ for all $i\in\mathbb{Z}_{\geq0}$, $y\in M(x)$. A mapping \tcb{$M$} is locally bounded (LB) if, for each bounded set $K$, $M(K):=\bigcup_{x\in K}M(x)$ is bounded.  Given a set $\mathcal{X} \subset \mathbb{R}^m$, the mapping $M$ is OSC and LB relative to $\mathcal{X}$ if the set-valued mapping from $\mathbb{R}^m$ to $\mathbb{R}^n$ defined by $M$ for $x \in \mathcal{X}$,  and by $\varnothing$ for $x \notin \mathcal{X}$, is OSC and LB at each $x \in \mathcal{X}$. The graph of $M$ is the set $\text{graph}(M) := \{(x, y) \in  \mathbb{R}^m \times \mathbb{R}^n: y \in M(x)\}$. Given a measurable space $(\Omega, \mathcal{F})$, a set-valued map $M: \Omega \rightrightarrows \mathbb{R}^n$ is said to be \textit{$\mathcal{F}$-measurable}, if for each open set $\mathcal{O} \subset \mathbb{R}^n$ we have $M^{-1}(\mathcal{O}) :=
	\{\omega\in\Omega:M(\omega)\cap\mathcal O\neq\varnothing\}
\in\mathcal F$.  
	
	\textbf{Stochastic Hybrid Dynamical Systems:} In this paper, we consider SHDS of the form
	\begin{subequations}\label{SHDS1}
		\begin{align}
			&y\in C,~~~~~~~~~\dot{y}\in F(y),\label{SHDS_flows0}\\
			&y\in D,~~~~~~y^+\in G(y,v^+),~~~v\sim \mu(\cdot),~~\label{SHDS_jumps0}
		\end{align}
	\end{subequations}
	where $F:\mathbb{R}^n\rightrightarrows\mathbb{R}^n$ is called the flow map, $G:\mathbb{R}^n\times\mathbb{R}^m\rightrightarrows\mathbb{R}^n$ is called the jump map, $C$ is the flow set, $D$ is the jump set, and $v^+$ is a place holder for a sequence $\{\bf{v}_k\}_{k=1}^{\infty}$ of independent, identically distributed {\em i.i.d.} random variables $\bf{v}_k:\Omega\to\mathbb{R}^m$, $k\in\mathbb{N}$, defined on a probability space $(\Omega,\mathcal{F},\mathbb{P})$. Thus, ${\bf v_k}^{-1}\tcb{(B)}:=\{\omega\in\Omega: \mathbf{v}_k(\omega)\in \tcb{B}\}\in\mathcal{F}$ for all $\tcb{B}\in\mathbf{B}(\tcb{\mathbb{R}^m})$, and $\mu:\mathbf{B}(\mathbb{R}^m)\to[0,1]$ is defined as $\mu(\tcb{B}):=\mathbb{P}\{\omega\in\Omega:{\bf v}_k(\omega)\in \tcb{B}\}$. 
	
	Following \cite{Teel_ANRC}, we denote a random solution
to \eqref{SHDS1} by $\mathbf y(\omega)$, where $\omega\in\Omega$.
Let $\mathcal F_0=\{\varnothing,\Omega\}$ and
$\mathcal F_\ell:=\sigma(\mathbf v_1,\ldots,\mathbf v_\ell)$
for $\ell\geq1$. The sequence $(\mathcal F_\ell)_{\ell\geq0}$
is the natural filtration of $(\mathbf v_\ell)_{\ell\geq1}$.  A stochastic hybrid arc is a mapping ${\bf y}$ from $\Omega$ to the set of hybrid arcs \cite[Ch. 2]{bookHDS}, such that the set-valued mapping from $\Omega$ to $\mathbb{R}^{n+2}$, given by  $\omega\mapsto \text{graph}({\bf y}(\omega)):=\big\{(t,j,z):\tilde{y}={\bf y}(\omega), (t,j)\in\text{dom}(\tilde{y}),z=\tilde{y}(t,j)\big\}$, is $\mathcal{F}$-measurable with closed values. Let $\text{graph}({\bf y}(\omega))_{\leq \ell}:=\text{graph}({\bf y} (\omega))\cap (\mathbb{R}_{\geq0}\times\{0,1,\ldots,\ell\}\times\mathbb{R}^n)$. An $\{\mathcal{F}_\ell\}_{\ell=0}^{\infty}$ adapted stochastic hybrid arc is a stochastic hybrid arc ${\bf y}$ such that the mapping $\omega\mapsto \text{graph}({\bf y}(\omega))_{\leq \ell}$ is $\mathcal{F}_\ell$ measurable for each $\ell \in\mathbb{Z}_{\geq0}$. An adapted stochastic hybrid arc ${\bf y}(\omega)$, or simply $\mathbf{y}_\omega$,  is a solution to \eqref{SHDS1} starting from $y_0$ denoted ${\bf y}_\omega\in \mathcal{S}_r(y_0)$ if: (1) $\mathbf{y}_\omega(0,0)=y_0$; (2) if $(t_1,j),(t_2,j)\in\text{dom}(\mathbf{y})$ with $t_1<t_2$, then for almost all $t\in[t_1,t_2]$, $\mathbf{y}_{\omega}(t,j)\in C$ and $\dot{\mathbf{y}}_\omega(t,j)\in F(\mathbf{y}_\omega(t,j))$; (3) if $(t,j),(t,j+1)\in\text{dom}(\mathbf{y}_\omega)$, then $\mathbf{y}_\omega(t,j)\in D$ and $\mathbf{y}_\omega(t,j+1)\in G(\mathbf{y}_\omega(t,j),\mathbf{v}_{j+1}(\omega))$. A random solution ${\bf y}_\omega$ is: a) almost surely {\em non-trivial} if its hybrid time domain contains at least two points almost surely; and b) almost surely  {\em  complete} if for almost every sample path $\omega\in \Omega$ the hybrid arc ${\bf y_\omega}$ has an unbounded time domain. 

	\vspace{-0.2cm} 
	\section{A Class of Singularly Perturbed SHDS}
	\label{section3}

    \vspace{-0.2cm}
	In this paper, we consider SHDS of the form
	\begin{subequations}\label{SPSHDS1}
		\begin{align}
			&(x,z)\in C,~~\dot{x}\in F_x(x,z),~~~\varepsilon\dot{z}\in F_z(x,z)
			\label{flow_dynamics1}\\
			&(x,z)\in D,~~(x^+,z^+)\in G(x,z,v),~~v\sim \mu,\label{jump_dynamics1}
		\end{align}
	\end{subequations}
	where $\varepsilon\in\mathbb{R}_{>0}$ is a small parameter, $x\in\mathbb{R}^{n_1}$ is the ``slow'' state, $z\in\mathbb{R}^{n_2}$ is the ``fast'' state, $F_x:\mathbb{R}^{n_1}\times\mathbb{R}^{n_2}\rightrightarrows\mathbb{R}^{n_1}$, $F_z:\mathbb{R}^{n_1}\times\mathbb{R}^{n_2}\rightrightarrows\mathbb{R}^{n_2}$ are set-valued mappings, and $C,D\subset\mathbb{R}^{n_1}\times\mathbb{R}^{n_2}$ define the flow set and the jump sets, respectively. To simplify our presentation, we will consider flow and jump sets of the form $C:=C_x\times C_z$ and $D:=D_x\times D_z$, where $C_x,D_x\subset\mathbb{R}^{n_1}$ and $C_z,D_z\subset\mathbb{R}^{n_2}$ are all closed sets. We will also use the notation
	$$F_{\varepsilon}(x,z):=F_x(x,z)\times \varepsilon^{-1}F_z(x,z),$$
	to denote the overall flow map, $y=(x,z)$ to denote the overall state, and $\mathcal V=\{\mathbf v_k(\omega):\omega\in\Omega,\ k\ge1\}$ to denote the set of all possible values that $v$ can take.
	
	\vspace{-0.2cm} 
	In system~\eqref{SPSHDS1}, \(v\) denotes a sequence of i.i.d.\ random variables with probability measure~\(\mu\). Since setting \(\varepsilon=0\) induces a singularity in the flow dynamics~\eqref{flow_dynamics1}, system~\eqref{SPSHDS1} constitutes a \emph{singularly perturbed stochastic hybrid dynamical system} (SP-SHDS). This general model encompasses several classes previously studied in the literature. For example, when \(G\) is independent of \(v\), it reduces to the deterministic hybrid systems analyzed in~\cite{SP2020Teel,RajebSP}. If, in addition, \(F_x\) and \(F_z\) are single-valued and continuous and the flow and jump sets for \(z\) are compact, the model specializes to those in~\cite{Sanfelice:11,PoTe17Auto}. When \(D=\emptyset\), it recovers singularly perturbed differential inclusions~\cite{HistorySingularPerturbations}, and if \(F_x\) and \(F_z\) are single-valued and locally Lipschitz, it reduces to the singularly perturbed ODEs studied in~\cite{Saberi}. Finally, for linear dynamics,~\eqref{SPSHDS1} recovers the classical linear singular perturbation models of~\cite[Ch.~1--6]{Kokotovic_SP_Book}.
	
	\vspace{-0.2cm} 
	The inclusion of stochastic dynamics in singularly perturbed hybrid inclusions introduces a rich range of behaviors that remain largely unexplored in the singular perturbation literature. In this work, we focus on \emph{causal} behaviors consistent with the notion of random solutions adopted here. \tcb{We also impose the following regularity conditions on the system data ~\((C,F_{\varepsilon},G,D,\mu)\)}.
	
    
	\vspace{-0.1cm}
\begin{definition}\label{definitionbasic1}
		The SP-SHDS \eqref{SPSHDS1} is said to satisfy the \textsl{Basic Conditions} if the following conditions hold: (a) $C$ and $D$ are closed,  $C\subset \text{dom}(F_{\varepsilon})$. (b) $F_{\varepsilon}$ is OSC and LB relative to $C$, and convex-valued on $C$, for any $\varepsilon>0$. (c) $G$ is LB relative to $D\times\mathbb{R}^m$ and the mapping $v\mapsto \text{graph}(G(\cdot,\cdot,v)):=\{(s_0,s_1,s_2)\in D\times\mathbb{R}^{n_1+n_2}:s_2\in G(s_0,s_1,v)\}$ is measurable with closed values, \tcb{and $G(x,z,v)\neq\emptyset$ for all $(x,z)\in D,\ v\in\mathbb R^m$.} \QEDB
	\end{definition}

    \vspace{-0.1cm} 
	In this paper, we make the standing assumption that all the SHDS satisfy the Basic Conditions. 
	
	%
	
	%

    \vspace{-0.1cm} 
	Following the classic paradigm of singular perturbations \cite{Kokotovic_SP_Book,Saberi,Teel96SP}, our goal is to inform the stability properties of \eqref{SPSHDS1} based on the properties of a simplified \emph{reduced} SHDS characterized by the ``steady-state'' behavior of the fast dynamics in \eqref{flow_dynamics1}, which describe the \emph{boundary-layer} dynamics of the system.

	\vspace{-0.1cm} 
	\subsection{Boundary Layer Dynamics}

	\vspace{-0.2cm} 
	In the proposed framework, the boundary layer dynamics of \eqref{SPSHDS1} ignore the jumps, and are given by
	\begin{equation}\label{BLD}
		(x,z)\in C,~~\dot{z}\in F_{z}(x,z),~~\dot{x}=0.
	\end{equation}
	Note that in \eqref{BLD}, the state $x$ is fixed. For these dynamics, we will assume the existence of a quasi steady-state map $\mathcal{M}$ that satisfies the following regularity assumption:

    \vspace{-0.2cm}
	\begin{assumption}\label{manifold}
		There exists an OSC and LB set-valued mapping $\mathcal{M}:\mathbb{R}^{n_1}\rightrightarrows\mathbb{R}^{n_2}$ such that $\mathcal{M}(x)\neq\emptyset$ and $\mathcal{M}(x)\subset C_z$ for all $x\in C_x$. Moreover, \tcb{there exists $\alpha\in\mathcal{K}$ such that} 

        \vspace{-0.6cm}\begin{equation}\label{preservingmap2}
			\min_{r_2\in\mathcal{M}(x)}\|r_2^*-r_2\|\leq \alpha(\|x-r_1^*\|),
		\end{equation}

        \vspace{-0.4cm}\noindent 
		for all $x\in\tcb{\mathbb{R}^{n_1}}$ and for all $r^*\in \{r^*=(r_1^*,r_2^*)\in\mathbb{R}^{n_1+n_2}:r^*_1\in\mathcal{A},r^*_2\in \mathcal{M}(r^*_1)\}$.
	\end{assumption}

	\vspace{-0.1cm}
    \tcb{For each fixed $x$, the quasi-steady-state map $\mathcal M(x)$ represents a target set for the boundary-layer dynamics. In the classical single-valued setting, this set is typically an equilibrium $m(x)$ satisfying $f_z(x,m(x))=0$. When $\mathcal M(x)=\{m(x)\}$ with $m$ continuous, the first part of Assumption 1 recovers the regularity requirements typically imposed on the quasi-steady-state map in the ODE literature \cite{Saberi,Kokotovic_SP_Book,TeelMoreauNesic2003}, while condition (4) reduces to the “stability-preserving” condition $|m(x)-m(x^*)|\leq\alpha(|x-x^*|)$ considered in \cite[Ch. 11, pp. 450]{KhalilBook}. For differential inclusions, the corresponding equilibrium condition is $0\in F_z(x,m(x))$. Here, we allow $\mathcal M(x)$ to be set-valued and do not require it to consist exclusively of equilibria. Instead, its regularity and associated stability properties are specified in the following assumptions}:

	\vspace{-0.1cm} 
	\begin{assumption}\label{Assumption2}
		There exists a locally Lipschitz and regular function $W:\mathbb{R}^{n_1}\times\mathbb{R}^{n_2}\to\mathbb{R}_{\geq0}$, functions $\alpha_1,\alpha_2\in\mathcal{G}_{\infty}$, a continuous function $\varphi_z:\mathbb{R}_{\geq0}\to\mathbb{R}_{\geq0}$, and $k_z>0$ such that:

        \vspace{-0.2cm}
		\begin{enumerate}[(a)]
			\item For all $y\in C\cup D\cup G(D\times\mathcal{V})$, 
			\begin{equation}\label{Kinftyboundsfast}
				\alpha_1(|z|_{\mathcal{M}(x)})\leq W(x,z)\leq \alpha_2(|z|_{\mathcal{M}(x)}).
			\end{equation}
			\item For all $y\in C$ and all $f_z\in F_{z}(x,z)$,
			\begin{equation}\label{boundaryflowsLyapunov}
				\max_{\substack{v\in \partial_z W(x,z)}} \langle v,f_z\rangle  \leq -k_z \varphi_z^2(|z|_{\mathcal{M}(x)}).
			\end{equation}
		\end{enumerate}
	\end{assumption}
	For each fixed $x$, and by using the results of \cite{Teel_ANRC}, the conditions of Assumption \ref{Assumption2} essentially establish uniform Lagrange stability of the compact set $\mathcal{M}(x)$ for the DI \eqref{BLD}. If, additionally, $\alpha_1,\alpha_2\in\mathcal{K}_{\infty}$ and $\varphi_z$ is positive definite, they imply uniform global asymptotic stability.
	
	\vspace{-0.1cm} 
	In some applications, the stochastic jumps of the fast state $z$ might also converge towards \tcb{$\mathcal{M}(x)$} in a probabilistic sense. In that case, the following Assumption will be considered.

	\vspace{-0.1cm} \begin{assumption}\label{assumptionjumpsfast}
		There exists $c_z>0$ and $\rho_z\in\mathcal{P}_s\mathcal{D}$ such that the function $W$ of Assumption \ref{Assumption2} satisfies:
		\begin{equation}
			\int_{\mathbb{R}^m}\sup_{\substack{g\in G(x,z,v)}}~W(g)\mu(dv)\leq W(x,z)-c_z\rho_z(|z|_{\mathcal{M}(x)}).
		\end{equation}
		for all $y\in D$. 
	\end{assumption}

	\vspace{-0.1cm} 
	In words, Assumption \ref{assumptionjumpsfast} asks that the worst case value of the function $W$ evaluated during jumps of the \emph{original} SHDS \eqref{SPSHDS1} does not increase in expectation.
	
	\vspace{-0.1cm} 
	\subsection{Reduced Stochastic Hybrid Dynamics}

	\vspace{-0.2cm} 
	Using Assumption \ref{manifold}, we define the \emph{reduced} SHDS:
	\begin{subequations}\label{SPSHDSreduced}
		\begin{align}
			&x\in C_x,~~~~\dot{x}\in \tilde{F}(x),\label{flow_dynamics1r}\\
			&x\in D_x,~~x^+\in \tilde{G}(x,v),~~v\sim \mu,\label{jump_dynamics1r}
		\end{align}
	\end{subequations}
	where $\tilde{F}$ and $\tilde{G}$ are the set-valued mappings:
	\begin{subequations}
		\begin{align}
			&\tilde{F}(x):=\overline{\text{co}}\left\{\tilde{f}_x\in\mathbb{R}^{n_1}:\tilde{f}_x\in F_x(x,z),z\in \mathcal{M}(x)\right\},\\
			&\tilde{G}(x,v):=\left\{s\in\mathbb{R}^{n_1}:(s,l)\in G(x,z,v),z\in D_z\right\}.
		\end{align}
	\end{subequations}
	%
	As mentioned before, we assume that system \eqref{SPSHDSreduced} satisfies the Basic Conditions. \tcb{Under Assumption~\ref{manifold}, the Basic Conditions
for \eqref{SPSHDS1} imply those for \eqref{SPSHDSreduced}
if $D_z$ is nonempty and either compact or $G$ is
independent of $z$.}

	Throughout, we consider a compact set $\mathcal{A}\subset\mathbb{R}^{n_1}$, an open and bounded set $\mathcal{O}_x\subset\mathbb{R}^{n_1}$, and a function $\varpi_x:\mathbb{R}^{n_1}\rightarrow\mathbb{R}_{\geq 0}$. In most cases, $\varpi_x(x)=|x|_{\mathcal{A}}$. To capture the dynamic properties of \eqref{SPSHDSreduced}, we will first consider the following assumption on the flows.
	
	\begin{assumption}\label{Assumption3}
		There exists a locally Lipschitz and regular function $V:\mathbb{R}^{n_1}\to\mathbb{R}_{\geq0}$, a $C^0$ function \tcb{$\varphi_x:\mathbb{R}^{n_1}\to\mathbb{R}_{\geq0}$}, $\alpha_3,\alpha_4\in \mathcal{G}_{\infty}$, and $k_x>0,\nu\geq0$, such that:

        \vspace{-0.2cm}
		\begin{enumerate}[(a)]
			\item For all $x\in C_x\cup D_x \cup \tilde{G}(D_x\times\mathcal{V})$, 
			\begin{equation}\label{slowKinfty}
				\alpha_3(|x|_{\mathcal{A}})\leq V(x)\leq \alpha_4(\varpi_x(x)).
			\end{equation}
			\item For all $x\in C_x$ and all $\tilde{f}_x\in \tilde{F}(x)$,
			\begin{equation}\label{reducedflowsLyapunov}
				\max_{\substack{v\in \partial V(x)}}\langle v,\tilde{f}_x\rangle \leq -k_x \varphi_x^2(x)+\nu \mathbb{I}_{\mathcal{O}_x}(x).
			\end{equation}
		\end{enumerate}
	\end{assumption}
	As in Assumption \ref{Assumption2}, the two conditions of Assumption \ref{Assumption3} essentially guarantee uniform Lagrange stability of the set $\mathcal{A}$ relative to $\varpi_x$ for the DI \eqref{flow_dynamics1r} whenever $\nu=0$, via \cite[Theorem 4.2]{Teel_ANRC}.  When $\nu>0$ and $\mathcal{A}=\overline{\mathcal{O}_x}$, condition \eqref{reducedflowsLyapunov} allows the study of recurrence properties for the open set $\mathcal{O}_x$. 
    
    In some applications, the reduced SHDS \eqref{SPSHDSreduced} might also exhibit suitable stability properties during the jumps. In that case, the following assumption will be used in conjunction with Assumption \ref{Assumption3}.
	\begin{assumption}\label{assumptionjumpsreduced}
		There exists a $C^0$ function $\rho_x:\mathbb{R}^{n_1}\to\mathbb{R}_{\geq0}$ and $c_x>0$ such that  
\begin{equation}\label{jumpsreducedstochastic}
			\int_{\mathbb{R}^m}\sup_{\substack{\tilde{g}\in \tilde{G}(x,v)}}~V(\tilde{g})\mu(dv)\leq V(x)-c_x\rho_x(x),
		\end{equation}
		for all $x\in D_x$. 
	\end{assumption}
	If Assumptions \ref{Assumption3} and \ref{assumptionjumpsreduced} hold with $\mathcal{A}=\overline{\mathcal{O}_x}$, $\varphi_x(x)^2\geq \tilde{c}_x>0$ for all $x\in C_x\setminus\mathcal{O}_x$, and $\rho_x(x)>0$ for all $x\in D_x\setminus\mathcal{O}_x$, 
	then we can directly conclude uniform global recurrence (see Definition~\ref{definition3}) of the set $\mathcal{O}_x$ relative to $\varpi_x$ for the \emph{reduced} SHDS \eqref{SPSHDSreduced} via \cite[Thm. 4.4]{Teel_ANRC}. Separately, if Assumptions~\ref{Assumption3} and~\ref{assumptionjumpsreduced} hold with $\alpha_3,\alpha_4\in\mathcal{K}_{\infty}$, $\rho_x,\varphi_x(\cdot)^2\in\mathcal{P}\mathcal{D}(\mathcal{A})$, $\nu=0$, and $\varpi_x$ is lower semicontinuous with sub-level sets that satisfy certain natural conditions \cite[Eq. (21)]{Teel_ANRC}, then we can conclude UGASp (see Definition~\ref{SHDS2def}) of the set $\mathcal{A}$ relative to $\varpi_x$ for the \emph{reduced} SHDS via \cite[Thm. 4.5]{Teel_ANRC}. 
	
	\vspace{-0.2cm} 
	Whether or not these stability properties will be preserved, in some sense, by the original SP-SHDS \eqref{SPSHDS1}, will depend on the time-scale separation, parameterized by the value of $\varepsilon$, as well as some additional ``interconnection'' conditions.
	\subsection{Interconnection Conditions}

    \vspace{-0.1cm}
	We consider two different types of interconnection conditions: one related to the flows \eqref{flow_dynamics1}, and the other one related to the jumps \eqref{jump_dynamics1}. 

    \vspace{-0.1cm}
	Below, in Assumptions \ref{assumption6}-\ref{interconnection2} the functions $V$, $W$, $\varphi_x$, and $\varphi_z$ are the same from Assumptions \ref{Assumption2} and \ref{Assumption3}.

    \vspace{-0.1cm}
	\begin{assumption}\label{assumption6}
		There exist constants $k_1,k_3>0$ and \tcb{$k_2,k_4\geq0$} such that:
		\begin{enumerate}[(a)]
			\item For all $y\in C$, and all $f_x\in F_x(x,z)$:
			\begin{equation}\label{coupled1Lyapu1}
				\begin{aligned}
					\max_{\substack{v\in \partial_x W(x,z)}} &\langle v, f_x \rangle \leq k_1\varphi_z(|z|_{\mathcal{M}(x)})\varphi_x(x)\\
					&+k_2\varphi_z^2(|z|_{\mathcal{M}(x)}) + k_4\varphi_z(|z|_{\mathcal{M}(x)}).
				\end{aligned}
			\end{equation}
			\item For all $y\in C$, and for all $f_x\in F_x(x,z)$, there exists $\tilde{f}_x\in \tilde{F}(x)$ such that
			\begin{equation}\label{coupled1Lyapu2}
				\begin{aligned}
					\max_{\substack{v\in \partial V(x)}} \langle v,f_x -\tilde{f}_x \rangle &\leq k_3 \varphi_z(|z|_{\mathcal{M}(x)})\varphi_x(x) \\
					&+ k_4\varphi_z(|z|_{\mathcal{M}(x)}). 
				\end{aligned}
			\end{equation}
		\end{enumerate}
		\vspace{-0.3cm}
	\end{assumption}
    When $D=G=\emptyset$, $\nu=k_4=0$, $\mathcal{A}=\{0\}$, $\varphi_x,\varphi_z\in\mathcal{P}\mathcal{D}$, and all the vector fields are single-valued and locally Lipschitz, the conditions of Assumptions \ref{Assumption2}, \ref{Assumption3}, and \ref{assumption6} essentially recover the conditions in \cite{Saberi} for asymptotic stability in singularly perturbed smooth ODEs. 

    \vspace{-0.1cm}
    The next two assumptions concern the behavior of the system during jumps and can be invoked separately.

    \vspace{-0.1cm}
\begin{assumption}\label{jumpcoupledcondition}
		There exist  \tcb{$k_5\geq0$} and a continuous function $\rho_5:\mathbb{R}^{n_1}\to\mathbb{R}_{\geq0}$, such that

        \vspace{-0.9cm}
        \begin{equation}\label{jumpcondition_expectation}
			\int_{\mathbb{R}^m}\sup_{\substack{g\in G(x,z,v)}}~W(g)\mu(dv)\leq W(x,z)+k_5\rho_5(x), 
		\end{equation}
		for all $y\in D$. 
	\end{assumption}
	    \vspace{-0.1cm}\begin{assumption}\label{interconnection2}
		There exist  \tcb{$k_6\geq0$} and a continuous function $\rho_6:\mathbb{R}_{\geq0}\to\mathbb{R}_{\geq0}$, such that for all $y\in D$:

        \vspace{-0.9cm}
		\begin{equation}		\int_{\mathbb{R}^m}\sup_{\substack{\tilde{g}_x\in \tilde{G}(x,v)}}~V(\tilde{g}_x)\mu(dv)\leq V(x)+k_6\rho_6(|z|_{\mathcal{M}(x)}).
		\end{equation} 
		\vspace{-0.2cm}
	\end{assumption}

	%
	%
	%
	%

    \vspace{-0.5cm}
	\subsection{Stability Results}
	\label{sec_results}

    \vspace{-0.3cm}
	To study the stability properties of the SP-SHDS \eqref{SPSHDS1}, we consider a class of composite regular, locally Lipschitz certificate functions of the form
	\begin{equation}\label{foster_function}
		E_{\theta}(y):=(1-\theta) V(x)+\theta W(x,z), 
	\end{equation}
	where $\theta\in(0,1)$, and $V,W$ come from Sections \ref{section3}.1-\ref{section3}.3. For convenience, we also introduce the set $L_{E_{\theta}}(c):=\{y\in\text{dom}~E_{\theta}:E_{\theta}(y)=c\}$, and the constant 

    \vspace{-0.8cm}
	\begin{align}\label{thetacond}
		\theta^*:= \frac{k_3}{k_1+k_3}.
	\end{align}

    \vspace{-0.4cm}
	\subsubsection{ \tcb{Asymptotic} Stability in Probability}
    \vspace{-0.2cm}
	Given a compact set $\mathcal{A}\subset\mathbb{R}^{n_1}$, we are interested in establishing suitable stochastic stability properties for the SP-SHDS \eqref{SPSHDS1} with respect to the set

    \vspace{-0.5cm}\begin{equation}\label{compactsetcomplete}
		\tilde{\mathcal{A}}:=\left\{y\in\mathbb{R}^{n_1+n_2}:x\in\mathcal{A},~z\in \mathcal{M}(x)\right\},
	\end{equation}

    \vspace{-0.2cm}
	which is assumed nonempty and is compact by compactness of $\mathcal A$ and the OSC and local boundedness of $\mathcal M$. In particular, we focus on the property of \emph{Uniform Global Stability in Probability}, introduced in \cite[Sec. 2.3]{Teel_ANRC}.
	
	\vspace{-0.1cm}
	\begin{definition}\label{SHDS2def}
		Consider the SP-SHDS \eqref{SPSHDS1}, the set $\mathcal{\tilde{A}}\subset\mathbb{R}^{n_1+n_2}$, and a function $\varpi:C\cup D\rightarrow\mathbb{R}_{\geq 0}$.

        \vspace{-0.3cm}
		\textbf{(C1)} The set $\mathcal{\tilde{A}}$ is said to be \textsl{Uniformly Lyapunov stable in Probability} relative to $\varpi$ if for each $\epsilon>0$ and $\rho>0$ there exists a $\delta>0$ such that every maximal random solution with $\varpi(\mathbf{y}_\omega(0,0))\leq \delta$ satisfies:
			\begin{align}\label{stabilityprobability}
				\mathbb{P}\big(\mathbf{y}_{\omega}(t,j)&\in\mathcal{\tilde{A}}+\epsilon\mathbb{B}^\circ, \notag\\
				&~\forall~ (t,j)\in\text{dom}(\mathbf{y}_{\omega})\big)\geq 1-\rho.
			\end{align}

        \vspace{-0.4cm}
         \textbf{(C2)} The set $\mathcal{\tilde{A}}$ is said to be \textsl{Uniformly Lagrange stable in probability} relative to $\varpi$ if for each $\delta>0$ and $\rho >0$, there exists $\epsilon>0$ such that inequality \eqref{stabilityprobability} holds for every maximal random solution with $\varpi(\mathbf{y}_\omega(0,0))\leq \delta$. 

	      \vspace{-0.2cm}
		\textbf{(C3)} The set $\mathcal{\tilde{A}}$ is said to be \textsl{Uniformly Globally Attractive in Probability} relative to $\varpi$ if for each $\epsilon>0,$ $\rho>0$, and $R>0$, there exists $T\geq 0$ such that for all random solutions  ${\bf y}_\omega$ satisfying $\varpi(\mathbf{y}_\omega(0,0))\leq R$: 
            \vspace{-0.3cm}
			\begin{align*}
				\mathbb{P}\Big(&\mathbf{y}_{\omega}(t,j)\in\mathcal{\tilde{A}}+\epsilon\mathbb{B}^\circ, \notag\\
				&~~~~~~~\forall~t+j\geq T,(t,j)\in \text{dom}(\mathbf{y}_{\omega})\Big)\geq 1-\rho.
			\end{align*}

        \vspace{-0.4cm}
		If conditions \textbf{(C1)}, \textbf{(C2)}, and \textbf{(C3)} hold, system \eqref{SPSHDS1} is said to render the set $\mathcal{\tilde{A}}$ \textsl{Uniformly Globally Asymptotically Stable in Probability} (UGASp) relative to $\varpi$. If $\varpi=|y|_{\tilde{\mathcal{A}}}$, we omit ``relative to $\varpi$" from the definition. 
	\end{definition}

    \vspace{-0.1cm}
	Our first result provides different sufficient conditions to guarantee UGASp of the set $\tilde{\mathcal{A}}$ for the SP-SHDS \eqref{SPSHDS1}.

    \vspace{-0.1cm}
	\begin{theorem}\label{theorem1}
		Let $\mathcal{A}\subset\mathbb{R}^{n_1}$ be compact, and $\tilde{\mathcal{A}}$ be given by \eqref{compactsetcomplete}. Suppose that Assumption \ref{manifold} holds, and that:

        \vspace{-0.3cm}
		\begin{enumerate}[(a)]
			\item Assumptions \ref{Assumption2}, \ref{Assumption3} and \ref{assumption6} hold with $\alpha_i\in\mathcal{K}_{\infty}$ for all $i\in\{1,2,3,4\}$, $\varphi_x\in\mathcal{P}\mathcal{D}(\mathcal{A})$, $\varphi_z\in\mathcal{P}\mathcal{D}$, $\varpi_x(x)=|x|_{\mathcal{A}}$, and $\nu=k_4=0$.
			\item There exists a function $\hat{\rho}\in\mathcal{P}\mathcal{D}(\tilde{\mathcal{A}})$ such that

            \vspace{-0.5cm}\begin{equation}\label{jump_conditionproof1}
	\int_{\mathbb{R}^m}\sup_{\substack{g\in G(x,z,v)}}~E_{\theta^*}(g)\mu(dv)\leq E_{\theta^*}(y)-\hat{\rho}(y),
			\end{equation}
			for all $y\in D$, where $E_{\theta^*}$ is given by \eqref{foster_function}.
			%
			%
		\end{enumerate}
		Then, there exists $\varepsilon^*>0$ such that, for all $\varepsilon\in(0,\varepsilon^*)$, the set $\tilde{\mathcal{A}}$ is UGASp for the \tcb{SP-SHDS \eqref{SPSHDS1}}. \\ \phantom{.}\hfill
	\end{theorem}

    \vspace{-0.3cm}
	\textbf{Proof:} We show that, under the \tcb{given} conditions, $E_{\theta^*}$ is a strong regular Lyapunov-Foster function for the SP-SHDS. First, by \cite[Lemma 9]{SPSHDS}, there exist  functions $\hat{\alpha}_1,\hat{\alpha}_2\in\mathcal{K}_{\infty}$ such that for all $y\in C\cup D\cup G(D\times\mathcal{V})$:
	\begin{equation}\label{kinftyproof1}
		\hat{\alpha}_1(|y|_{\tilde{\mathcal{A}}})\leq E_{\theta^*}(y)\leq\hat{\alpha}_2(|y|_{\tilde{\mathcal{A}}}),
	\end{equation}
	and by \cite[Prop. 3.2.6]{Clarke}, $E_{\theta}$ is regular for each $\theta\in(0,1)$. \tcb{Since $\partial E_{\theta}(y)
=(1-\theta)\bigl(\partial V(x)\times\{0_{n_2}\}\bigr)
+\theta\partial W(x,z)$}, and $\partial W(x,\tcb{z})\subset \partial_x W(x,\tcb{z})\times \tcb{\partial_z W(x,z)}$, using \eqref{foster_function} it follows that for all $y\in C$ and all $f=(f_x,\varepsilon^{-1}f_z)$ with $f_x\in F_x(x,z)$ and \tcb{$f_z\in F_z(x,z)$}, we have that:
	\begin{align*}
		\max_{\substack{e\in \partial E_{\theta}(y)}}\langle e,f\rangle\leq &(1-\theta) \max_{\substack{v\in \partial V(x)}}\langle v,f_x\rangle+\theta \max_{\substack{v\in \partial_x W(y)}}\langle v,f_x\rangle\\
		&+\frac{1}{\varepsilon}\theta \max_{\substack{v\in \partial_z W(y)}}\langle v,f_z\rangle.
	\end{align*}
	Adding and subtracting terms, we obtain:

    \vspace{-1cm}
	\begin{align*}
		\max_{\substack{e\in \partial E_{\theta}(y)}}&\langle e,f\rangle\leq\theta \max_{\substack{v\in \partial_x W(y)}}\langle v,f_x\rangle +\frac{\theta}{\varepsilon} \max_{\substack{v\in \partial_z W(y)}}\langle v,f_z\rangle\\
		&+(1-\theta)\left(\max_{\substack{v\in \partial V(x)}}\langle v,f_x \rangle-\max_{\substack{v\in \partial V(x)}}\langle v,\tilde{f}_x \rangle\right)\\
		&+(1-\theta) \max_{\substack{v\in \partial V(x)}}\langle v,\tilde{f}_x\rangle,
	\end{align*}
	for all $f\in F_{\varepsilon}(x,z)$, where $\tilde{f}_x$ comes from \eqref{coupled1Lyapu2}.  Using \eqref{boundaryflowsLyapunov}, \eqref{reducedflowsLyapunov}, \eqref{coupled1Lyapu1}, \eqref{coupled1Lyapu2}, the subadditive property of the function $S_{\partial V(x)}(a):=\max_{v\in\partial V(x)}\langle v,a\rangle$, and combining terms:
	\begin{align*}
		\max_{\substack{e\in \partial E_{\theta}(y)}}\langle e,f\rangle \leq & -(1-\theta)k_x\varphi_x^2(x)\\
		&+(\theta k_1+(1-\theta)k_3)\varphi_z(|z|_{\mathcal{M}(x)})\varphi_x(x)\\
		&-\theta\left(\frac{k_z}{\varepsilon}-k_2\right)\varphi_z^2(|z|_{\mathcal{M}(x)}),
	\end{align*}
	for all $f\in F_{\varepsilon}(x,z)$ (recall that $\nu=k_4=0$). 
	Using similar computations such as those in \cite{Saberi}, 
	it follows that, provided  $\varepsilon\in(0,\varepsilon^*)$ with \tcb{$\varepsilon^*=\frac{k_x k_z}{2(k_2k_x+k_1k_3)}$}, there exists $\lambda>0$ such that for all $y\in C$, and all $f\in F_{\varepsilon}(x,z)$:
\begin{equation}\label{flowproof1}
    \max_{e\in\partial E_{\theta^*}(y)}
    \langle e,f\rangle
    \leq -\lambda p(y),~
    p(y):=\varphi_x(x)^2+
    \varphi_z(|z|_{\mathcal M(x)})^2,
\end{equation}
for all $y=(x,z)\in C$ and all $f\in F_\varepsilon(y)$. \tcb{By \cite[Lemma 10]{SPSHDS}, there exists a continuous positive definite function
$\rho_F:\mathbb R_{\geq0}\to\mathbb R_{\geq0}$ such that
$\rho_F(|y|_{\tilde{\mathcal A}})\leq\lambda p(y)$ for all
$y\in C$}. Moreover, by item~(b), we have that
\begin{equation}\label{jumpproof1}
    \int_{\mathbb R^m}
        \sup_{g\in G(x,z,v)}E_{\theta^*}(g)\,\mu(dv)
        -E_{\theta^*}(y)
    \leq -\hat\rho(y)
\end{equation}
for all $y=(x,z)\in D$. Define $\rho(y):=
    \min\left\{
        \rho_F(|y|_{\tilde{\mathcal A}}),\,
        \hat\rho(y)
    \right\}$. The function $\rho$ is continuous and positive definite
with respect to $\tilde{\mathcal A}$. By
\eqref{flowproof1} and \eqref{jumpproof1}, both the flow
derivative and the expected jump increment of
$E_{\theta^*}$ are bounded above by $-\rho(y)$. Together with
\eqref{kinftyproof1}, these inequalities show that
$E_{\theta^*}$ is a Lyapunov--Foster function with respect to
$\tilde{\mathcal A}$. The result follows from
\cite[Thm.~4.5]{Teel_ANRC}.
\hfill$\blacksquare$ 
	
	\vspace{-0.1cm}
	The result of Theorem \ref{theorem1} relies on showing that $E_{\theta^*}$ is a
	\emph{strong regular Lyapunov-Foster function} for the SP-SHDS \eqref{SPSHDS1} whenever $\varepsilon$ is sufficiently small. However, in some applications it might not be easy to find functions $V,W$ that satisfy all the conditions of Theorem \ref{theorem1}. The next result relaxes some
of these conditions:
	%
	
	\begin{theorem}\label{theorem2}
		Let $\mathcal{A}\subset\mathbb{R}^{n_1}$ be compact, and $\tilde{\mathcal{A}}$ be given by \eqref{compactsetcomplete}. Suppose that Assumption \ref{manifold} holds, $\theta^*$ is given by \eqref{thetacond}, and
		\begin{enumerate}[(a)]
			\item Assumptions \ref{Assumption2}, \ref{Assumption3}, and \ref{assumption6} hold with $\alpha_i\in\mathcal{K}_{\infty}$, for all $i\in\{1,2,3,4\}$, $\varphi_x\in\mathcal{P}_s\mathcal{D}(\mathcal{A})$, $\varphi_z\in\mathcal{P}_s\mathcal{D}$, $\varpi_x(x)=|x|_{\mathcal{A}}$, and $\nu=k_4=0$.
			%
			%
			\item At least one of the following conditions holds:
			\begin{enumerate}[(1)]
				\item \hspace{-0.1cm}Assumptions \ref{assumptionjumpsreduced} and \ref{jumpcoupledcondition} with $\nu=0$, $\rho_x=\rho_5$, and
\begin{equation}\label{relaxcondition1}
					\frac{k_3k_5}{k_1}<c_x.
				\end{equation}
				\item Assumptions \ref{assumptionjumpsfast} and \ref{interconnection2} with $\nu=0$, $\rho_z=\rho_6$, and
\begin{equation}\label{relaxcondition2}
					\frac{k_1k_6}{k_3}<c_z.
				\end{equation}
			\end{enumerate}
			\item For all sufficiently small $\varepsilon>0$ and every $c>0$, there is no almost surely complete random solution of \eqref{SPSHDS1} that remains in $L_{E_{\theta^*}}(c)$ almost surely.
		\end{enumerate}
		Then, there exists $\varepsilon^*>0$ such that, for all $\varepsilon\in(0,\varepsilon^*)$, the system \eqref{SPSHDS1} renders the set $\tilde{\mathcal{A}}$ UGASp.
	\end{theorem}

    \vspace{-0.1cm}
	\textbf{Proof:} By \cite[Lemma 9]{SPSHDS}, $E_{\theta^*}$ satisfies
the comparison bounds in \eqref{kinftyproof1}.
However, the quadratic-form estimate \eqref{flowproof1} now gives  $\max_{e\in\partial E_{\theta^*}(y)}\langle e,f\rangle
    \leq
    -\lambda\bigl(
        \varphi_x(x)^2+
        \varphi_z(|z|_{\mathcal M(x)})^2
    \bigr)
    \leq 0$, for all $y=(x,z)\in C$,  $f\in F_\varepsilon(y)$, and all
$\varepsilon\in(0,\varepsilon^*)$. During jumps, we have that for all $y\in D$:
	\begin{align*}
&\int_{\mathbb{R}^m}\sup_{\substack{g\in G(x,z,v)}}~E_{\theta^*}(g)\mu (dv)\\
&=\int_{\mathbb{R}^m}\sup_{\substack{(g_x,g_z)\in G(x,z,v)}}~\Big((1-\theta^*)V(g_x)+\theta^* W(g_x,g_z)\Big)\mu (dv)
\end{align*}
\begin{align*}
        &\leq  \int_{\mathbb{R}^m} \Bigg(\sup_{\substack{(g_x,g_z)\in G(x,z,v)}}(1-\theta^*)V(g_x)\\
&~~~~~~~~~~~~~~~~~~+\sup_{\substack{(g_x,g_z)\in G(x,z,v)}}\theta^* W(g_x,g_z)\Bigg)\mu (dv)\\
&= \int_{\mathbb{R}^m} \sup_{\substack{(g_x,g_z)\in G(x,z,v)}}(1-\theta^*)V(g_x)\mu (dv)\\
        &~~~~~~~~~~~~~~~~~~+\int_{\mathbb{R}^m} \sup_{\substack{(g_x,g_z)\in G(x,z,v)}}\theta^* W(g_x,g_z)\mu (dv)\\
		&\tcb{\leq} (1-\theta^*)\int_{\mathbb{R}^m} \sup_{\substack{\tilde{g}\in \tilde{G}(x,v)}}V(\tilde{g})\mu (dv)\\
		&~~~~~~~~~~~~~~~~~~+\theta^*\int_{\mathbb{R}^m} \sup_{\substack{(g_x,g_z)\in G(x,z,v)}} W(g_x,g_z)\mu (dv).
	\end{align*}
	where we used the monotonicity and linearity properties of the expectation, and the definition of system \eqref{SPSHDSreduced}. \tcb{In particular, the last inequality follows because, for every $(x,z)\in D$
and every $v\in\mathbb R^m$, we have $\{g_x:(g_x,g_z)\in G(x,z,v)
           \text{ for some }g_z\}
    \subseteq \widetilde G(x,v)$}. Using Assumptions \ref{assumptionjumpsreduced} and \ref{jumpcoupledcondition}, we get that for all $y\in D$:

    \vspace{-0.9cm}
	\begin{align*}
		\int_{\mathbb{R}^m}\sup_{\substack{g\in G(x,z,v)}}~E_{\theta^*}(g)\mu (dv)&\leq (1-\theta^*) \big(V(x)-c_x\rho_x(x)\big)\\
		&~~~+\theta^* W(x,z)+\theta^* k_5\rho_5(x)\\
		&=E_{\theta^*}(y)-c\rho_x(x)\leq E_{\theta^*}(y),
	\end{align*}
	where we used $c:=(1-\theta^*)c_x-\theta^* k_5$, the assumption that $\rho_x=\rho_5$, the definition of $\theta^*$, and the fact that $c>0$ due to \eqref{relaxcondition1}. Thus, $E_{\theta^*}$ satisfies
	\begin{align*}
		&\max_{\substack{e\in \partial E_{\theta^*}(y)}}\langle e,f\rangle \leq 0,~~\forall~y\in C,~\forall~f\in F_{\varepsilon}(x,z),\\
		&\int_{\mathbb{R}^m}\sup_{\substack{g\in G(x,z,v)}}\hspace{-0.3cm}E_{\theta^*}(g)\mu (dv)-E_{\theta^*}(y)\leq0,~~\forall~y\in D.
	\end{align*}
	For the case when Assumptions \ref{assumptionjumpsfast} and \ref{interconnection2} hold with $\rho_z=\rho_6$, we get that for all $y\in D$:
	\begin{align*}
		\int_{\mathbb{R}^m}\sup_{\substack{g\in G(x,z,v)}}&\hspace{-0.3cm}E_{\theta^*}(g)\mu (dv)\leq \theta^*\big(W(x,z)-c_z\rho_z(|z|_{\mathcal{M}(x)})\big)\\
		&+(1-\theta^*) \big(V(x)+k_6\rho_6(|z|_{\mathcal{M}(x)})\big)\\
		&=E_{\theta^*}(y)-\tilde{c}\rho_6(|z|_{\mathcal{M}(x)})\leq E_{\theta^*}(y),
	\end{align*}
	where $\tilde{c}:=\theta^* c_z-(1-\theta^*)k_6$, which is positive due to \eqref{relaxcondition2}. Since there does not exist an almost surely complete random solution $\bf y$ that remains in a non-zero level set of $E_{\theta^*}$ almost surely, UGASp of $\tilde{\mathcal{A}}$ follows directly via the stochastic hybrid invariance principle \cite[Thm. 8]{AnanathTACRecurrence}. \hfill $\blacksquare$

    \vspace{0.0cm}
	\subsubsection{Uniform Global Recurrence}

	\vspace{-0.1cm}
	The second main property that we study in this paper is the property of uniform global recurrence, introduced in \cite[Sec. 2.4]{Teel_ANRC}. This property is commonly studied in stochastic systems for which a stable set might not exist. 
	\vspace{-0.1cm}
	\begin{definition}\label{definition3}
		An open, bounded set $\mathcal{O}\subset\mathbb{R}^{n_1+n_2}$ is \emph{uniformly globally recurrent} (UGR) relative to $\varpi$ for the SP-SHDS \eqref{SPSHDS1} if there are no finite escape times for \eqref{flow_dynamics1} and for each $\rho>0$ and $R>0$ there exists $\tau\geq0$ such that for every maximal random solution with $\varpi({\bf{y}}_{\omega}(0,0))\leq R$,
	\begin{align*}
\mathbb P\Big(
&\big(\operatorname{graph}(\mathbf y_\omega)
      \subseteq\Gamma_{<\tau}\times\mathbb R^{n_1+n_2}\big)\\
&{}\lor
 \big(\operatorname{graph}(\mathbf y_\omega)
      \cap(\Gamma_{\leq\tau}\times\mathcal O)
      \neq\varnothing\big)
\Big)\geq 1-\rho,
\end{align*}
where $\operatorname{graph}(\mathbf y_\omega)
=\{(t,j,\mathbf y_\omega(t,j)):
(t,j)\in\operatorname{dom}(\mathbf y_\omega)\}$ and
$\Gamma_{<\tau}:=\{(t,j):t+j<\tau\}$, with
$\Gamma_{\leq\tau}$ defined analogously. \\ \phantom{.}\hfill
	\end{definition}

    \vspace{-0.3cm}
	Loosely speaking, the property of recurrence rules out finite escape times along flows and requires that, for each $\rho\in(0,1)$, solutions stop
or hit $\mathcal O$ with probability at least $1-\rho$ by a hybrid-time
bound that is uniform over each initial $\varpi$-sublevel set. This property is studied mainly for open sets in order to guarantee suitable uniformity and robustness properties. In particular, we study the recurrence of the set
	\begin{align}\label{eq:recurr_set_331}
		\mathcal{O}_{\chi}=\{(x,z)\in \mathbb{R}^{n_1+n_2} ~|~ x\in\mathcal{O}_x,~|z|_{\mathcal{M}(x)}<\chi\},
	\end{align}
	where $\mathcal{O}_x$ is an open and bounded set, and $\chi > 0$. In this way, the set $\mathcal{O}_{\chi}$ will be open and bounded by construction. 

    \vspace{-0.1cm}
	\begin{theorem}\label{theorem3}
		Let $\mathcal{O}_x\subset\mathbb{R}^{n_1}$ be an open and bounded set and suppose the following conditions hold:

        \vspace{-0.1cm}
		\begin{enumerate}[(a)]
			\item Assumptions \ref{manifold}, \ref{Assumption2}, \ref{Assumption3} and \ref{assumption6} hold with $\varphi_z\in\mathcal{PD}$, $\varphi_x(x)\geq\tilde{c}_x>0$ for all $x\in C_x$, and $\mathcal{A}=\overline{\mathcal{O}_x}$.
			\item There exist a continuous function
$\hat{\rho}:\mathbb{R}^{n_1+n_2}\to\mathbb{R}_{>0}$
and constants $\nu,\chi>0$, with $\nu$ at least as large as the
flow bias in Assumption~\ref{Assumption3}$\text{(b)}$, such that
$\tcb{m_\chi:=\inf_{s\geq\chi}\varphi_z(s)>0}$ and

            \vspace{-0.5cm}
			\begin{equation*}
				\int_{\mathbb{R}^m}\sup_{\substack{g\in G(x,z,v)}}~E_{\theta^*}(g)\mu(dv)\leq E_{\theta^*}(y)-\hat{\rho}(y)+\nu\mathbb{I}_{\mathcal{O}_{\chi}}(y),
			\end{equation*}
			for all $y\in D$, where $E_{\theta^*}$ is given by \eqref{foster_function} and \eqref{thetacond}.
			%
			%
		\end{enumerate}
		Then, there exists $\varepsilon^*>0$ such that, for all
$\varepsilon\in(0,\varepsilon^*)$, the set $\mathcal O_\chi$ is UGR
relative to $\varpi(y)=\|y\|$ for the SP-SHDS \eqref{SPSHDS1}.
	\end{theorem}
    The following example illustrates Theorem \ref{theorem3}.
    \begin{example}\normalfont
        Let $x,z\in\mathbb{R}$, $v\in\{-1,1\}$, $C_x=\mathbb{R}$, $D_x=\mathbb{Z}_{\geq 0}$, $C_z=D_z=\mathbb{R}$, and consider the SP-SHDS \eqref{SPSHDS1} with
        \begin{align*}
            F_x(x,z)&=\{z\}, \quad F_z(x,z)=\{-(z+x)\}, \quad v\sim \mu,\\
            G(x,z,v)&=\{\max\{0,x+v\}\}\times\{-\max\{0,x+v\}\},
        \end{align*}
        where $\mu(\{-1\})=\frac{17}{20}$ and $\mu(\{1\})=\frac{3}{20}$. The associated reduced-order SHDS is \eqref{SPSHDSreduced} with $\tilde{F}(x)=\{-x
        \}$ and $\tilde{G}(x,v)=\{\max\{0,x+v\}\}$. Let $\mathcal O_x=(-1,1)$, $\mathcal A=[-1,1]$, and
$\chi=1$. Define $\mathcal M(x)=\{-x\}$ and $V(x)=\tfrac12 x^2$, $W(x,z)=\tfrac12(z+x)^2$. Assumption~\ref{manifold} holds with $\alpha(r)=r$. The comparison bounds on $V$ and $W$ hold with
$\alpha_i(r)=r^2/2$ for $i=1,2,3,4$ and
$\varpi_x(x)=|x|$, since $|x|_{\mathcal A}\leq |x|$
and $|z|_{\mathcal M(x)}=|z+x|$. Choose $\varphi_x(x)=\sqrt{x^2+\tfrac1{10}}$, $\varphi_z(s)=s$, $k_x=\tfrac12$, $k_z=1$, $\nu=\tfrac1{10}$. In particular, $\varphi_x(x)\geq 1/\sqrt{10}$ for all $x$,
and $m_\chi=\inf_{s\geq1}\varphi_z(s)=1$. Along the
reduced and boundary-layer flows, respectively, we have
\begin{align*}
    \nabla V(x)(-x)
    &=-x^2
    \leq -\tfrac12\varphi_x(x)^2
           +\tfrac1{10}\mathbb I_{\mathcal O_x}(x),\\
    \partial_z W(x,z)\bigl(-(z+x)\bigr)
    &=-(z+x)^2.
\end{align*}
The first inequality follows by considering $|x|<1$
and $|x|\geq1$ separately.  For $f_x=z$ and $\tilde f_x=-x$, we also have $\partial_x W(x,z)f_x\leq |z+x||x|+(z+x)^2$ and
    $\nabla V(x)(f_x-\tilde f_x)\leq |x||z+x|$. Since $|x|\leq\varphi_x(x)$, Assumption~\ref{assumption6}
holds with $k_1=k_2=k_3=1$ and $k_4=0$. Thus
$\theta^*=1/2$, and the composite function is $E_{\theta^*}(x,z)=\tfrac14\bigl[x^2+(z+x)^2\bigr]$. At a jump, $x^+=\max\{x+v,0\}$ and $z^+=-x^+$.
Therefore, $E_{\theta^*}(x^+,z^+)=\tfrac14(x^+)^2$,
and, for every $y=(x,z)\in D$, we have $\int_{\mathbb R}\sup_{g\in G(x,z,v)}
        E_{\theta^*}(g)\,\mu(dv)
    =E_{\theta^*}(y)-\tilde\rho(y)$, where $\tilde\rho(x,z)
    =\tfrac14\bigl[x^2+(z+x)^2\bigr]-\tfrac3{80}\max\{x+1,0\}^2
        -\tfrac{17}{80}\max\{x-1,0\}^2$. \tcb{To bound this decrement, let
$\tilde\rho(x,z)=\tfrac14(z+x)^2+q(x)$, where
\[
q(x)=
\begin{cases}
    x^2/4, & x\leq -1,\\
    (17x^2-6x-3)/80, & -1<x<1,\\
    7x/20-1/4, & x\geq1.
\end{cases}
\]}
The middle quadratic has minimum $-3/68$ at $x=3/17$,
and $q(x)\geq 1/10$ for $|x|\geq1$. If $|x|<1$ and
$(x,z)\notin\mathcal O_\chi$, then $|z+x|\geq1$, so
$\tilde\rho(x,z)\geq 1/4-3/68=7/34>1/10$.
Consequently, $\tilde\rho\geq1/10$ outside
$\mathcal O_\chi$ and $\tilde\rho\geq-3/68$ globally.
Since $-3/68\geq-1/20$, the jump condition in
Theorem~\ref{theorem3}(b) holds with
$\hat\rho(y)=1/20$ and $\nu=1/10$. It follows from Theorem~\ref{theorem3} that
$\mathcal O_\chi$ is UGR for all sufficiently small
$\varepsilon>0$. We note, however, that no compact set is UGASp for this system. Indeed, similar to the example in \cite[p. 306]{subbaraman2015robustness}, the system admits almost surely complete \emph{discrete} random solutions that are unbounded.         \QEDB 
    \end{example}

    \vspace{-0.2cm}
	As in Theorem \ref{theorem2}, we can relax the decrease conditions on $E_{\theta^*}$ to establish UGR of $\mathcal{O}_{\chi}$ for sufficiently small $\varepsilon$, provided that certain solutions of \eqref{SPSHDS1} can be ruled out. 

	\vspace{0.0cm}
	\begin{theorem}\label{theorem4}
		Let $\mathcal{O}_x\subset\mathbb{R}^{n_1}$ be open and
bounded, and let $\chi>0$. Suppose that \tcb{$m_\chi:=\inf_{s\geq\chi}\varphi_z(s)>0$}, and that the following conditions hold: 

        \vspace{-0.2cm}
		\begin{enumerate}[(a)]
			\item Assumptions \ref{manifold}, \ref{Assumption2}, \ref{Assumption3} and \ref{assumption6} hold with $\nu>0$, $\varphi_z\in\mathcal{PD}$, $\varphi_x(x)\geq\tilde{c}_x>0$ for all $x\in C_x$, and $\mathcal{A}=\overline{\mathcal{O}_x}$.
			%
			%
			\item For all $y\in D\backslash\mathcal{O}_{\chi}$, we have that
			\begin{equation*}
				\int_{\mathbb{R}^m}\sup_{\substack{g\in G(x,z,v)\cap(\mathbb{R}^{n_1+n_2}\backslash\mathcal{O}_{\chi})}}~E_{\theta^*}(g)\mu(dv)\leq E_{\theta^*}(y),
			\end{equation*}
			where $E_{\theta^*}$ is given by \eqref{foster_function}. The supremum is taken to be zero when $G(x,z,v)\cap
(\mathbb R^{n_1+n_2}\setminus\mathcal O_\chi)=\varnothing.$
			\item For all sufficiently small $\varepsilon>0$ and every
$c\geq0$, there is no almost surely complete random
solution of \eqref{SPSHDS1} that remains in $L_{E_{\theta^*}}(c)
    \cap(\mathbb R^{n_1+n_2}\setminus\mathcal O_\chi)$ almost surely.
		\end{enumerate}
        \vspace{-0.2cm}
		Then, there exists $\varepsilon^*>0$ such that, for all $\varepsilon\in(0,\varepsilon^*)$, the SP-SHDS \eqref{SPSHDS1} renders the set $\mathcal{O}_{\chi}$ UGR relative to some $\varpi:C\cup D\rightarrow\mathbb{R}_{\geq 0}$.
	\end{theorem}
	\textbf{Proof of Theorem \ref{theorem3}:} \tcb{Fix $(x_0,z_0)\in\mathcal O_\chi$ and choose
$w_0\in\mathcal M(x_0)$ with $\|z_0-w_0\|<\chi$.
Since $x_0\in\mathcal O_x\subset\mathcal A$,
Assumption~\ref{manifold} gives $|z|_{\mathcal M(x)}
    \leq \|z-w_0\|+\alpha(\|x-x_0\|)$. The right-hand side is less than $\chi$ for all
$(x,z)$ sufficiently close to $(x_0,z_0)$. Since $\mathcal O_x$ is open, these points also satisfy
$x\in\mathcal O_x$. Thus $\mathcal O_\chi$ is open.
Its boundedness follows from boundedness of
$\mathcal O_x$ and local boundedness of $\mathcal M$.} Next, for any fixed $\theta\in(0,1)$, the lower bounds on
$V$ and $W$ imply that boundedness of $E_\theta(y)$
bounds both $|x|_{\mathcal A}$ and $|z|_{\mathcal M(x)}$
on the flow, jump, and reachable states. Compactness
of $\mathcal A$ and local boundedness of $\mathcal M$
then imply boundedness of $y$. These bounded sublevel
sets give a class-$\mathcal G_\infty$ lower comparison
bound in terms of $|y|_{\overline{\mathcal O_\chi}}$.
Continuity of $E_\theta$ gives a class-$\mathcal G_\infty$
upper comparison bound in terms of $\|y\|$, by taking
a continuous nondecreasing unbounded majorant of
$\max_{\|w\|\leq r}E_\theta(w)$.
	
    By Assumptions \ref{Assumption2}, \ref{Assumption3} and \ref{assumption6}, and following similar steps of the proof of Theorem \ref{theorem1}, we obtain that
		\begin{align*}
			\max_{\substack{e\in \partial E_{\theta}(y)}}\langle e,f&\rangle \leq    \nu \mathbb{I}_{\mathcal{O}_x}(x) -(1-\theta)k_x\varphi_x^2(x)\\
			&+(\theta k_1+(1-\theta)k_3)\varphi_z(|z|_{\mathcal{M}(x)})\varphi_x(x)\\
			&-\theta\left(\frac{k_z}{2\varepsilon}-k_2\right)\varphi_z(|z|_{\mathcal{M}(x)})^2 \\
			&+k_4 \varphi_z(|z|_{\mathcal{M}(x)})- \frac{k_z \theta }{2\varepsilon}\varphi_z(|z|_{\mathcal{M}(x)})^2,
		\end{align*}
        for any $\theta\in(0,1)$ and any $\varepsilon>0$.
Fix $\theta=\theta^*$ and define $\varepsilon_1:=\frac{k_xk_z}
    {4(k_xk_2+k_1k_3)}$. For $\varepsilon>0$, consider the matrix
\[
Q_\varepsilon:=
\begin{bmatrix}
(1-\theta^*)k_x
&
-\dfrac{\theta^*k_1+(1-\theta^*)k_3}{2}
\\[1ex]
-\dfrac{\theta^*k_1+(1-\theta^*)k_3}{2}
&
\theta^*\left(\dfrac{k_z}{2\varepsilon}-k_2\right)
\end{bmatrix}.
\]
\tcb{Using $\theta^*=k_3/(k_1+k_3)$, we obtain $\det Q_\varepsilon
=
\frac{k_1k_3}{(k_1+k_3)^2}
\left[
k_x\left(\frac{k_z}{2\varepsilon}-k_2\right)
-k_1k_3
\right]$. Consequently, $Q_{\varepsilon_1}$ is positive definite.
Moreover, $Q_\varepsilon-Q_{\varepsilon_1}$ is positive
semidefinite for every $\varepsilon\in(0,\varepsilon_1)$.
Thus, setting $\lambda:=\lambda_{\min}(Q_{\varepsilon_1})>0$, the preceding flow estimate gives
\begin{align*}
\max_{e\in\partial E_{\theta^*}(y)}
\langle e,f\rangle
\leq{}&
-\lambda\left(
\varphi_x(x)^2+
\varphi_z(|z|_{\mathcal M(x)})^2
\right)
+\nu\mathbb I_{\mathcal O_x}(x)
\\
&+k_4\varphi_z(|z|_{\mathcal M(x)})
-\frac{k_z\theta^*}{2\varepsilon}
\varphi_z(|z|_{\mathcal M(x)})^2,
\end{align*}
for all $y\in C$, all $f\in F_\varepsilon(y)$,
and all $\varepsilon\in(0,\varepsilon_1)$. Define $\tilde{\mathcal O}_\chi
:=
\left\{
(x,z)\in C\cup D:
x\in\mathcal O_x,\ 
|z|_{\mathcal M(x)}\geq\chi
\right\}$. For every $y=(x,z)\in C$, we have $\mathbb I_{\mathcal O_x}(x)
=
\mathbb I_{\mathcal O_\chi}(y)
+
\mathbb I_{\tilde{\mathcal O}_\chi}(y)$. Set $s:=|z|_{\mathcal M(x)}$ and define $R_F(y):=
-\frac{\lambda\tilde c_x^2}{2}
+\nu\mathbb I_{\tilde{\mathcal O}_\chi}(y)
+k_4\varphi_z(s)
-\frac{k_z\theta^*}{2\varepsilon}\varphi_z(s)^2$. Since $\varphi_x(x)\geq\tilde c_x$, the flow estimate
implies $\max_{e\in\partial E_{\theta^*}(y)}
\langle e,f\rangle
\leq
-\frac{\lambda}{2}\varphi_x(x)^2
-\lambda\varphi_z(s)^2
+\nu\mathbb I_{\mathcal O_\chi}(y)
+R_F(y)$. If $y\in\tilde{\mathcal O}_\chi$, then $s\geq\chi$
and hence $\varphi_z(s)\geq m_\chi>0$.
Define $\varepsilon_2
    :=
    \frac{\theta^*k_zm_\chi^2}
    {2(k_4m_\chi+\nu)}$. For $\varepsilon\in(0,\varepsilon_2)$, we have $\frac{\theta^*k_z}{2\varepsilon}
\varphi_z(s)^2-k_4\varphi_z(s)-\nu
\geq
\varphi_z(s)^2
\left(
\frac{\theta^*k_z}{2\varepsilon}
-\frac{k_4}{m_\chi}
-\frac{\nu}{m_\chi^2}
\right)
\geq0$. Therefore, $R_F(y)\leq-\frac{\lambda\tilde c_x^2}{2}$. If $y\notin\tilde{\mathcal O}_\chi$, completing the
square gives $k_4\varphi_z(s)
-\frac{\theta^*k_z}{2\varepsilon}\varphi_z(s)^2
\leq
\frac{\varepsilon k_4^2}{2\theta^*k_z}$. Define $\varepsilon_3:=\frac{\theta^*k_z\lambda\tilde c_x^2}{2k_4^2}$ if $k_4>0$, and $\varepsilon_3:=+\infty$ if $k_4=0$.
%
Then, we get $R_F(y)\leq-\frac{\lambda\tilde c_x^2}{4}$ for all $\varepsilon\in(0,\varepsilon_3)$. Combining both cases, let $\varepsilon^*
    :=\min\{\varepsilon_1,\varepsilon_2,\varepsilon_3\}$, and $\tilde\rho(y):=\frac{\lambda\tilde c_x^2}{4}$}, which is continuous, strictly positive, and satisfies
\[
\max_{e\in\partial E_{\theta^*}(y)}
\langle e,f\rangle
\leq
-\tilde\rho(y)+\nu\mathbb I_{\mathcal O_\chi}(y),
\]
for all $y\in C$, all $f\in F_\varepsilon(y)$,
and all $\varepsilon\in(0,\varepsilon^*)$. Similarly, during jumps we have that
	\begin{equation*}
		\int_{\mathbb{R}^m}\sup_{\substack{g\in G(x,z,v)}}~E_{\theta^*}(g)\mu(dv)-E_{\theta^*}(y)\leq -\hat{\rho}(y)+\nu\mathbb{I}_{\mathcal{O}_{\chi}}(y).
	\end{equation*}
	and since $\hat{\rho}$ is also continuous and positive, it follows that $E_{\theta^*}$ is a Lagrange-Foster function for the SP-SHDS \eqref{SPSHDS1} with respect to the open set $\mathcal{O}_{\chi}$. Since $\mathcal O_\chi$ is open and bounded, \cite[Thm.~4.4]{Teel_ANRC} gives UGR of $\mathcal O_\chi$ relative to $\varpi(y)=\|y\|$. \hfill$\blacksquare$

	\textbf{Proof of Theorem \ref{theorem4}:} Following the flow estimates in the proof of
Theorem~\ref{theorem3}, using
\tcb{$m_\chi=\inf_{s\geq\chi}\varphi_z(s)>0$},
there exists $\varepsilon^*>0$ such that $\max_{e\in\partial E_{\theta^*}(y)}
\langle e,f\rangle\leq0$, for all $y\in C\setminus\mathcal O_\chi$,
all $f\in F_\varepsilon(y)$, and all
$\varepsilon\in(0,\varepsilon^*)$. Also, as shown in the proof of Theorem~\ref{theorem3}, the set
$\mathcal O_\chi$ is open and bounded. \tcb{The lower bounds on $V$ and $W$, compactness of $\mathcal A$,
and local boundedness of $\mathcal M$ imply that $E_{\theta^*}$
has bounded sublevel sets on $C\cup D\cup G(D\times\mathcal V)$.
Continuity gives a class-$\mathcal G_\infty$ upper bound in
$\|y\|$. Since $E_{\theta^*}$ is nonincreasing during flows
outside $\mathcal O_\chi$, its value along each flow segment
is bounded by the larger of its initial value and its maximum
on $\overline{\mathcal O_\chi}$. Bounded sublevel sets therefore
exclude finite escape along flows. Choose $\varepsilon^*>0$ so that the exterior flow estimate
and item~(c) hold for all $\varepsilon\in(0,\varepsilon^*)$.
The flow estimate gives nonincrease of $E_{\theta^*}$ along flows
outside $\mathcal O_\chi$. For jumps from $D\setminus\mathcal O_\chi$,
item~(b) bounds the expected supremum of $E_{\theta^*}$ over
post-jump states outside $\mathcal O_\chi$ by its pre-jump value.
The supremum is taken to be zero when all post-jump states
lie in $\mathcal O_\chi$.
Item~(c) excludes almost surely complete random solutions
confined to any exterior level set. Hence,
\cite[Thm.~9]{AnanathTACRecurrence} yields UGR of
$\mathcal O_\chi$ relative to $\varpi(y)=\|y\|$.}
\hfill $\blacksquare$
	%
    
	%
    
    \tcb{The stability and recurrence results presented in Theorems \ref{theorem1}-\ref{theorem4} use the same type of composite auxiliary-function construction, but impose different requirements on its bounds and evolution. For asymptotic stability, the composite function satisfies class-$(\mathcal K_\infty)$ bounds relative to the target compact set and decrease conditions ensuring stability and convergence under sufficient time-scale separation. For recurrence, the comparison-function bounds are weaker, and the required decrease is imposed outside the target open set. In the next section, we present different applications of these results.}

	\section{Applications to Switching Systems and Feedback Optimization}
	\label{sec_apps}

   \vspace{-0.2cm}
We illustrate the preceding results with four applications.

    \vspace{-0.2cm}
	\subsection{Switching Linear Systems with Stochastic Spontaneous Mode Transitions}\label{subsec1}
	Let ${x}:=(\xi,q,\tau)\in\mathbb{R}^n\times\mathcal{Q}\times[0,T]$, where $\mathcal{Q}:=\{1,\dots,N\}$ and $T>0$. Let ${z}\in\mathbb{R}^{n_2}$, $v:=(v_1,v_2)\in\mathbb{R}^2$ and consider the SP-SHDS $\mathcal{H}_{\varepsilon}$ given by \eqref{SPSHDS1}, 
	where the maps $F_{x}$, $F_{{z}}$, and $G$ are 
	\begin{align*}
		F_{{x}}({x},{z})&:=\{(\tilde{A}_q \xi + B_q {z}, 0,s)~|~s\in[-\eta,0]\},\\
		F_{{z}}({x},{z})&:= \{H \xi + L {z} \},\\
		G({x},z,v)&:=\{(\xi,v_1,v_2)\}\times\{{z}\},
	\end{align*}

    \vspace{-0.6cm}
	for some $\eta\geq0$, matrices $\tilde{A}_q,B_q,H,L$ of appropriate dimensions, and where the flow set $C:=C_x\times C_z$ and the jump set $D:=D_x\times D_z$ are defined by
	\begin{align*}
		C_x&:=\mathbb{R}^n\times\mathcal{Q}\times[0,T], & C_z:&= \mathbb{R}^{n_2}\\
		D_x&:=\mathbb{R}^n\times\mathcal{Q}\times\{0\}, & D_z:&= \mathbb{R}^{n_2}.
	\end{align*}

    \vspace{-0.5cm}
	In this model, $v_1$ is a random variable, independent of $v_2$, and distributed according to a discrete probability measure with a finite support defined by
	$$\mu_1:=\sum_{q\in\mathcal{Q}}\lambda_q\delta_q, \quad \sum_{q\in\mathcal{Q}}\lambda_q = 1,$$

    \vspace{-0.5cm}
	where, for each $q\in\mathcal{Q}$, $\lambda_q\geq 0$ is the probability of mode $q$, and $\delta_q$ is the Dirac measure on $\mathbb{R}$ centered at the point $q$. 
	On the other hand, $v_2$ is a non-negative scalar random variable, independent from $v_1$, and distributed according to the uniform probability distribution, denoted $\mu_2$, on the compact connected interval $[0,T]$, i.e.,

    \vspace{-0.8cm}
	\begin{align*}
		\mu_2(\mathrm{d}v_2):=\begin{cases} T^{-1}\mathrm{d}v_2  \, & v_2\in[0,T] \\ 
			0 & \text{otherwise.} \end{cases}
	\end{align*}
	Thus, the random vector $v=(v_1,v_2)$ is distributed according to the product measure $\mu = \mu_1 \times \mu_2$. System $\mathcal{H}_{\varepsilon}$ models a class of singularly perturbed switching linear system with ``spontaneous" mode transitions. At each jump, the timer is reset to a value drawn uniformly from \([0,T]\), independently of previous draws. Note that the inter-jump times depend on the selected timer velocity. 

	\vspace{-0.2cm}
	Associated with $\mathcal{H}_{\varepsilon}$ is the reduced order SHDS $\mathcal{H}_0$ defined by \eqref{SPSHDSreduced},
	where $\tilde{F}$ and $\tilde{G}$ are
	\begin{align*}
		\tilde{F}({x})&:=\{(A_q \xi, 0,s)~|~s\in[-\eta,0]\}, \\
		\tilde{G}({x},v)&:=\{(\xi,v_1,v_2)\},
	\end{align*}
	and $A_q:= \tilde{A}_q - B_q L^{-1}H$. To ensure that the SHDS $\mathcal{H}_{0}$ is well-behaved, we impose the following assumption.

    \vspace{-0.2cm}\begin{assumption}\label{asmp:hurwitz_Aq_1}
		There exists $\sigma>0$, a family of symmetric positive definite matrices $\{P_q\}_{q\in\mathcal{Q}}$, and a positive definite matrix $P_{z}$, such that
		\begin{subequations}\label{eq:matrix_inequalities_1}
			\begin{align}
				A_q^\top P_q + P_q A_q + \sigma \eta T^{-1} P_q &\prec 0, \label{eq28a}\\
				\sigma^{-1}\log(1+\sigma)P- P_q &\prec 0,\label{eq28b}\\
				L^\top P_{z} + P_{z} L &\prec 0,\label{eq29c}
			\end{align}
		\end{subequations}
		for all $q\in\mathcal{Q}$, where $P=\sum_{q\in\mathcal{Q}}\lambda_q P_q$. 
	\end{assumption}

    \vspace{-0.2cm}
	\begin{remark}\label{remark1}\normalfont
		If $\eta$ is a free parameter, and the matrix $L$, as well as the family of matrices $\{A_q\}_{q\in\mathcal{Q}}$, are Hurwitz, then \eqref{eq:matrix_inequalities_1} is always feasible. Indeed, inequality \eqref{eq28b} is equivalent to  $P^{-\frac{1}{2}}P_q P^{-\frac{1}{2}} - \sigma^{-1} \log(1+\sigma)I \succ 0$. By computation $\lim_{\sigma\nearrow \infty} \sigma^{-1} \log(1+\sigma) = 0$,
		which implies that, for any $\epsilon>0$, there exists $\underline{\sigma}(\epsilon)>0$ such that

        \vspace{-0.7cm}
		\begin{align*}
			|\sigma^{-1} \log(1+\sigma)|<\epsilon, \quad \forall \sigma>\underline{\sigma}(\epsilon).
		\end{align*}

        \vspace{-0.4cm}
		However, since $P_q$ may always be chosen as positive definite, there exists $\epsilon>0$ such that $P^{-\frac{1}{2}}P_q P^{-\frac{1}{2}} - \epsilon I \succ 0$
		which implies that, for such an $\epsilon>0$, we have

        \vspace{-0.7cm}
		\begin{align*}
			P^{-\frac{1}{2}}P_q P^{-\frac{1}{2}} - \sigma^{-1} \log(1+\sigma)I &\succ 0, \quad \forall \sigma>\underline{\sigma}(\epsilon),
		\end{align*}

        \vspace{-0.4cm}
		and, therefore, inequality  \eqref{eq28b} is satisfied. On the other hand, since $A_q^\top P_q + P_q A_q \prec 0$, 
		then, for any $\sigma>0$ there exists $\overline{\eta}(\sigma)>0$ such that inequality \eqref{eq28a} is satisfied for all $\eta\in(0,\overline{\eta}(\sigma))$. \hfill $\Box$
	\end{remark}

    \begin{figure}[t!]
			    \centering
			    \includegraphics[width=0.95\linewidth]{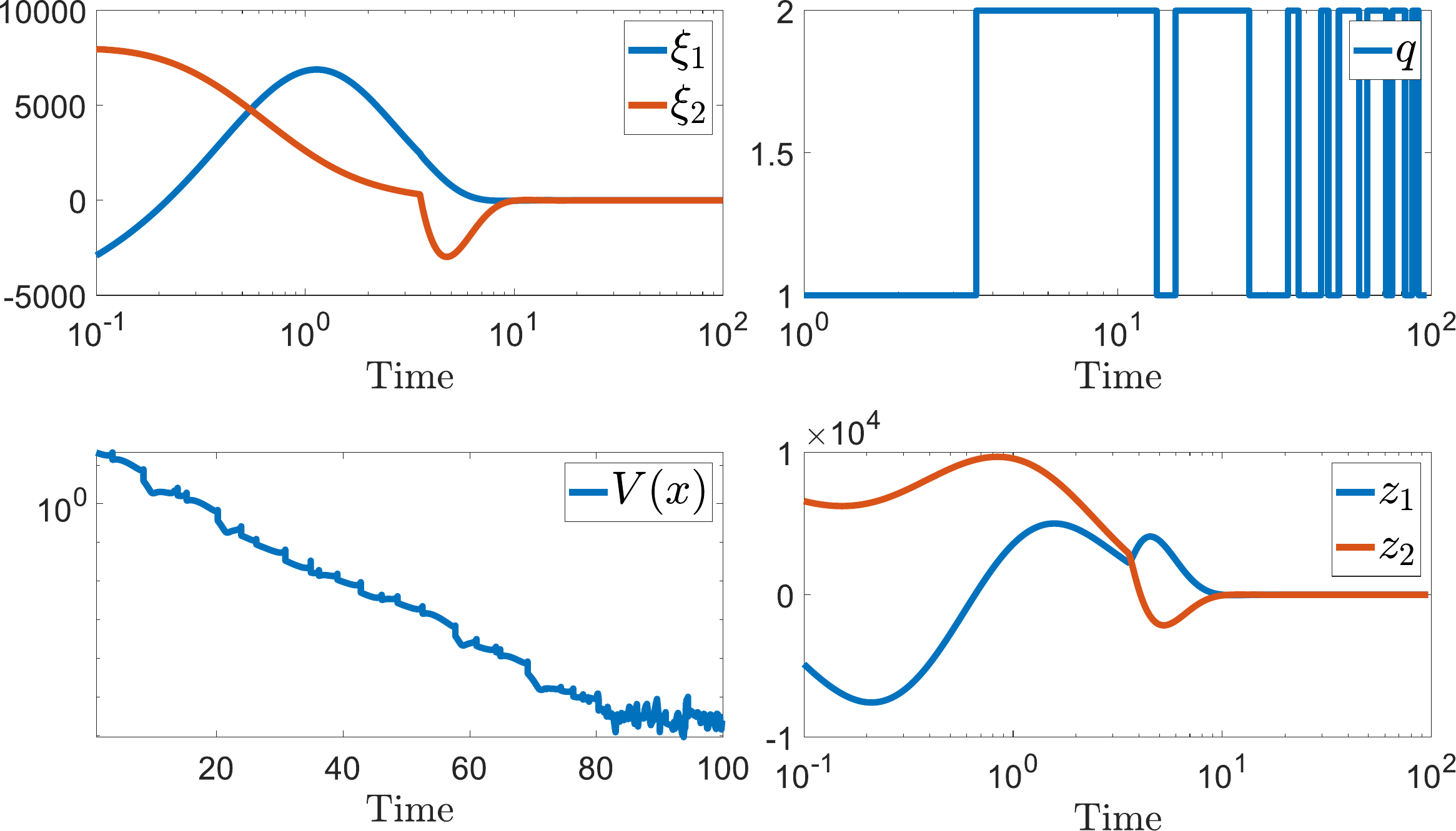}
			    \caption{Simulation results for switching system in Example \ref{example2}.}
			    \label{fig:exmp_1}
			\end{figure}
            
	Under Assumption \ref{asmp:hurwitz_Aq_1}, we have the following result.
	\begin{proposition}\label{thm:switching_systems}
		Let Assumption \ref{asmp:hurwitz_Aq_1} be satisfied and define the compact set $\tilde{\mathcal{A}}:= \mathcal{A}\times \{0\}$, where \tcb{$\mathcal{A}=\{0\}\times\mathcal{Q}\times[0,T]\subset C_x\cup D_x$}. Then, there exists $\varepsilon^*\in(0,\infty)$ such that, for all $\varepsilon\in (0,\varepsilon^*)$, the SP-SHDS $\mathcal{H}_\varepsilon$ renders UGASp the set $\tilde{\mathcal{A}}$. 
	\end{proposition}
The following example \tcb{illustrates} Proposition 1.
    \begin{example}\normalfont\label{example2}
		Let $\xi\in\mathbb{R}^2$, $z\in\mathbb{R}^2$, $\mathcal{Q}=\{1,2\}$, $T=2$, and $\tilde{A}_1=[-2,2;-1,0]$, $B_1=[0,1;1,0]$, $\tilde{A}_2=[-2,-1;-2,-2]$, $B_2=[0,1;-1,0]$, $L=[-1,0;0,-1]$, $H=[1,-1;1,1]$, where we used the so-called Matlab notation.
		Direct computation gives $A_1=[-1,3;0,-1]$, $A_2=[-1,0;-3,-1]$,
		%
		%
        both of which are Hurwitz. The pair admits no common quadratic Lyapunov function with strict decrease. We consider the uniform distribution over $\mathcal{Q}$, i.e., $\lambda_1=\lambda_2=\frac{1}{2}.$
		For $\sigma=18$, direct computation shows that $P_1=[0.50, 0.75;0.75,2.75]$, $P_2=[2.75,-0.75;-0.75,0.50]$
		%
		satisfy the second inequality in \eqref{eq:matrix_inequalities_1}. Moreover, we have that $A_1^\top P_1 + P_1 A_1= A_2^\top P_2 + P_2 A_2 = [-1, 0;0,-1]\prec 0$. Consequently, it can be shown that, for any $\eta < 0.035$, the first inequality in \eqref{eq:matrix_inequalities_1} is satisfied. Although this bound is conservative, it is worth mentioning that, unless $\eta=0$, the SP-SHDS $\mathcal{H}_{\varepsilon}$ admits sample paths with switching signals that have arbitrarily fast switching. However, this does not affect UGASp. Figure~\ref{fig:exmp_1} shows a sample path for $\varepsilon=0.1$. \hfill \QEDB
	\end{example}
    
	\textbf{Proof of Proposition 1:}
	Let $W:C\cup D\rightarrow \mathbb{R}_{\geq 0}$ be
	\begin{align*}
		W({x},{z}):=({z}+L^{-1}H\xi)^\top P_{{z}} ({z}+L^{-1}H\xi),
	\end{align*}
    \tcb{which is well-defined due to \eqref{eq29c}.} Let $V:C_x\cup D_x\rightarrow\mathbb{R}_{\geq 0}$ be given by
	\begin{align*}
		V({x}):=(\sigma\tau T^{-1}+1)^{-1} \xi^\top P_q \xi.
	\end{align*}
	Clearly, $W$ is positive definite and radially unbounded with respect to the slow manifold $\mathcal M(x):=\{-L^{-1}H\xi\}$. Indeed, since $|{z}|_{\mathcal{M}({x})} = \|{z} +L^{-1} H \xi \|$, we have that $\alpha_1(|{z}|_{\mathcal{M}({x})})\leq W({x},{z})\leq \alpha_2(|{z}|_{\mathcal{M}({x})})$, where the functions $\alpha_1$ and $\alpha_2$, given by $\alpha_1(r):=\lambda_{\min}(P_{z}) r^2$ and $\alpha_2(r):=\lambda_{\max}(P_{z}) r^2$, are both $\mathcal{K}_\infty$ functions. 
	In addition, by computation
	under Assumption \ref{asmp:hurwitz_Aq_1}, $W$ satisfies 
	\begin{align*}
		\max_{h\in \partial_{z} W({x},{z})}\langle h,f_{{z}}\rangle &\leq - k_{z} \varphi_{z}(|{z}|_{\mathcal{M}({x})})^2,
	\end{align*}
	for each $({x},{z})\in C$ and all $f_{{z}}\in F_{z}({x},{z})$, where $k_{z}$ and $\varphi_{z}$ are defined by $k_{z}:= -\lambda_{\max}(L^\top P_{{z}} + P_{{z}} L),~ \varphi_z(r):=|r|$, where, by virtue of Assumption \ref{asmp:hurwitz_Aq_1}, $k_{z}$ is strictly positive. 
	On the other hand, it is clear that $V$ is also positive definite and radially unbounded with respect to the compact set $\mathcal{A}$. Indeed, since $|{x}|_{\mathcal{A}} = \|\xi\|$, we have that $\alpha_3(|{x}|_{\mathcal{A}})\leq V({x}) \leq \alpha_4(|{x}|_{\mathcal{A}})$,
	%
	%
	for all ${x}\in C_x\cup D_x$, where $\alpha_3$ and $\alpha_4$ are defined by $\alpha_3(r):= (\sigma+1)^{-1}\min_{q\in\mathcal{Q}}\lambda_{\min}(P_q) r^2$,  $\alpha_4(r):=  \max_{q\in\mathcal{Q}}\lambda_{\max}(P_q) r^2$, which are both $\mathcal{K}_\infty$ functions.
	In addition, for each ${x}\in C_x$ and all $f_{x}\in \tilde{F}({x})$, we compute that
	\begin{align*}
		\max_{h\in \partial_{x} V({x})}&\langle h,f_{x}\rangle \leq (\sigma\tau T^{-1}+1)^{-1}\xi^\top(A_q^\top P_q + P_q A_q)\xi \\ 
		&\!\!\!+ \sigma \eta T^{-1}(\sigma\tau T^{-1}+1)^{-2} \xi^\top P_q \xi \\
		&\!\!\!\leq (\sigma\tau T^{-1}+1)^{-1} \xi^\top (A_q^\top P_q + P_q A_q + \sigma \eta T^{-1} P_q) \xi.
	\end{align*}
	Therefore, the function $V$ satisfies the inequality
	\begin{align*}
		\max_{h\in \partial_{x} V({x})}\langle h,f_{x}\rangle\leq - k_{x} \varphi_{x} ({x})^2 \leq 0,
	\end{align*}
	for all ${x}\in C_x$ and all $f_{x}\in \tilde{F}({x})$, with $\varphi_{x}({x}) := |{x}|_{\mathcal{A}}$, and
	\begin{align*}
		k_{x}&:= - (\sigma+1)^{-1} \max_{q\in\mathcal{Q}}\lambda_{\max}(A_q^\top P_q + P_q A_q + \sigma \eta T^{-1}P_q).
	\end{align*}
	Under Assumption \ref{asmp:hurwitz_Aq_1}, the constant $k_{x}$ is also strictly positive. 
	Moreover, direct computation gives that, for every ${x}\in D_x$, the function $V$ satisfies
	\begin{align*}
		\int_{\mathbb{R}^2}\sup_{g\in \tilde{G}({x},v)}V(g) &~\mu(\mathrm{d}v)-V({x}) \\
		&=  \xi^\top(\sigma^{-1}\log(1+\sigma)P - P_q) \xi,
	\end{align*}
	where $P=\sum_{q\in\mathcal{Q}}\lambda_{q}P_q$. Therefore, $V$ satisfies
	\begin{align*}
\int_{\mathbb{R}^2}\sup_{g\in\tilde{G}  (x,v)}V(g)\,\mu(\mathrm{d}v)
\leq V(x)-c_x\rho_x(x),
	\end{align*}
	for all ${x}\in D_x$, where $\rho_{x}({x}):= |{x}|_{\mathcal{A}}^2$, and the constant $c_{x}:= - \max_{q\in\mathcal{Q}}\lambda_{\max}(\sigma^{-1}\log(1+\sigma)P- P_q)$ is also strictly positive. For each $({x},{z})\in C$ and for all $f_{{x}}\in F_{x}({x},{z})$, we compute that
	\begin{align*}
		\max_{h\in \partial_{x} W({x},{z})} &\langle h,f_{x}\rangle \\
		&\!\!\!\!\!\!\!\!\!\!\!\!\!\!\!\!\!\!\!\!=2 \langle H^\top L^{-\top}P_{{z}} ({z}+L^{-1}H\xi), A_q \xi + B_q ({z}+L^{-1}H\xi)\rangle \\
		&\!\!\!\!\!\!\!\!\!\!\!\!\!\!\!\leq 2 \sigma_{\max}(A_q)\sigma_{\max}(H^\top L^{-\top}P_{{z}}) |{x}|_{\mathcal{A}} |{z}|_{\mathcal{M}({x})} \\
		&\!\!\!\!\!\!\!\!\!\!\!\!\!\!\!+2 \lambda_{\max}(\mathrm{Sym}(B_q^\top H^\top L^{-\top}P_{{z}})) |{z}|_{\mathcal{M}({x})}^2,
	\end{align*}
	and that, with $\tilde{f}_{{x}} = f_{{x}}+(A_q \xi-(\tilde{A}_q \xi + B_q {z}),0,0) \in \tilde{F}({x})$, 
	\begin{align*}
		\max_{h\in\partial_x V(x)}\langle h,f_x-\tilde f_x\rangle
&=\frac{\langle 2P_q\xi,B_q(z+L^{-1}H\xi)\rangle}{1+\sigma\tau/T}\\
		&\leq 2\sigma_{\max}(P_q) \sigma_{\max}(B_q) |{x}|_{\mathcal{A}} |{z}|_{\mathcal{M}({x})}.
	\end{align*}
	Hence, for each $({x},{z})\in C$ and every $f_{{x}}\in F_{x}({x},{z})$, there exists $\tilde{f}_{{x}} \in \tilde{F}({x})$ such that
	\begin{gather*}
		\max_{h\in \partial_{x} W({x},{z})} \langle h,f_{x}\rangle \leq k_1 \varphi_{x}({x}) \varphi_{z}(|{z}|_{\mathcal{M}({x})})+ k_2\varphi_{z}(|{z}|_{\mathcal{M}({x})})^2,\\
		\max_{h\in \partial_{x} V({x})} \langle h,f_{x} - \tilde{f}_{x}\rangle \leq k_3 \varphi_{x}({x})\varphi_{z}(|{z}|_{\mathcal{M}({x})}),
	\end{gather*}
	Choose $k_1,k_3>0$ and $k_2\geq0$ such that 
    \begin{align*}
    k_1&\geq 2\|H^\top L^{-\top}P_z\|\max_{q\in\mathcal Q}\|A_q\|\\
    k_2&\geq 2\max_{q\in\mathcal Q}\lambda_{\max}
       (\mathrm{Sym}(B_q^\top H^\top L^{-\top}P_z))\\
       k_3&\geq 2\max_{q\in\mathcal Q}\bigl(\|P_q\|\|B_q\|\bigr),
    \end{align*}
where $\|\cdot\|$ denotes the spectral norm. Finally, for all $({x},{z})\in D$, we compute the worst-case expectation $\int_{\mathbb{R}^2}\sup_{g\in G({x},{z},v)} W(g) ~\mu(\mathrm{d}v) =  W({x},{z})$
	%
	%
	which implies that, for all $({x},{z})\in D$,
	\begin{align*}
		\int_{\mathbb{R}^2}\sup_{g\in G({x},{z},v)} W(g) ~\mu(\mathrm{d}v) =  W({x},{z}) + k_5 \rho_5({x}),
	\end{align*}
	with $k_5:=0$ and any non-negative function $\rho_5({x})$. In particular, the same is true if we take $\rho_5({x}):=\rho_{x}({x})$. By combining all of the above, it follows that the SP-SHDS $\mathcal{H}_{\varepsilon}$ satisfies Assumptions \ref{manifold},\ref{Assumption2},\ref{Assumption3},\ref{assumptionjumpsreduced},\ref{assumption6}, and \ref{jumpcoupledcondition} in Section \ref{section3}. Therefore, using $E_{\theta^*}(x,z)= (1-\theta^*) V (x) + \theta^* W(x,z)$,
	%
	%
	a computation similar to the proof of Theorem 2 in Section \ref{sec_results} shows that there exists $\lambda > 0$ such that
	\begin{align*}
		\max_{h\in \partial E_{\theta^*}(x,z)} \langle h,f\rangle \leq - \lambda (\varphi_x(x)^2 + \varphi_z(|z|_{\mathcal{M}(x)})^2),
	\end{align*}
    for all $(x,z)\in C$, all $\varepsilon\in (0,\varepsilon^*)$ and all $f\in F_x(x,z)\times\varepsilon^{-1}F_z(x,z)$, and $$\int_{\mathbb{R}^2}\sup_{g\in G(x,z,v)}E_{\theta^*}(g)\,\mu(\mathrm{d}v)
\leq E_{\theta^*}(x,z),$$
	%
	%
	for all $(x,z)\in D$. To show that item (c) in Theorem \ref{theorem2} holds, we prove the following claim.

    \vspace{-0.2cm}\begin{claim}\label{claim:no_complete_solution}
		An almost surely complete solution that remains in a non-zero level set of $E_{\theta^*}$ almost surely does not exist.
	\end{claim}

    \vspace{-0.2cm}
    \textbf{Proof of Claim~\ref{claim:no_complete_solution}:}
\tcb{On a positive level set, $E_{\theta^*}$ strictly decreases
along every positive-duration flow segment. A complete
solution confined to such a level set must therefore jump
infinitely often without flowing. Almost surely, every timer
reset is positive, and another jump requires the timer to
flow to zero. This is a contradiction.} \hfill$\square$

	Thus item~(c) holds, and, by Theorem~\ref{theorem2}, the set $\tilde{\mathcal A}$ is UGASp. \hfill $\blacksquare$

    \vspace{-0.1cm}
	\subsection{Heavy-Ball Feedback Optimization with Spontaneous Momentum Resets}

    \vspace{-0.2cm}
	Let $T>0$, $x=(u,p,\tau)\in \mathbb{R}^n\times\mathbb{R}^n\times[0,T]$, $z\in\mathbb{R}^{n_2}$, and consider the SP-SHDS $\mathcal{H}_{\varepsilon}$ defined by \eqref{SPSHDS1}
	%
	where the maps $F_{x}$, $F_{{z}}$, and $G$ are given by
	\begin{align*}
		\!F_{{x}}({x},{z})&:=\{(p,-\beta p -\nabla \phi_u(u) - K \nabla\phi_y(L {z}+ d),1)\},\\
		F_{{z}}({x},{z})&:= \{A z + B u \},\\
		G(x,z,v)&:=\{(u,sp,v):s\in[0,\rho]\}\times\{z\},
	\end{align*}
	where $\beta>0$, $\rho\in[0,1)$, and $C=C_x\times\mathbb R^{n_2}$, $D=D_x\times\mathbb R^{n_2}$, with
	\begin{align*}
C_x&=\mathbb{R}^n\times\mathbb{R}^n\times[0,T], & D_x&=\mathbb{R}^n\times\mathbb{R}^n\times\{T\},
	\end{align*}
	and $v$ is a scalar random variable distributed according to some probability distribution $\mu$ supported on the compact interval $[0,T]$. The SP-SHDS $\mathcal{H}_{\varepsilon}$ models a feedback interconnection between a linear plant with state $z$ and input $u$, and a heavy-ball feedback optimization algorithm \cite{OnlinePovedaHybrid,FlorianSteadyState} that incorporates resets of the momentum $p$ via the jump rule $p^+ \in [0,\rho] p$ whenever the state of the timer $\tau$ satisfies $\tau = T$. The timer resets to the random draw $v$, while the momentum reset is set-valued. Related algorithms with stochastic restarting have been studied in the optimization literature \cite{gupta2022stochastic,belan2018restart,pokutta2020restarting,wang2022scheduled,fu2025hamiltonian}. We make the following regularity assumption:
\begin{assumption}\label{asmp:stable_plant}
		There exists $P\succ 0$ such that $A^\top P + P A \prec 0$, and there exists $m_{\phi},L_{\phi}$,$L_{\phi_y}>0$, such that
		\begin{align*}
			\phi(u)&:=\phi_u(u) + \phi_y(H u + d), & H&:=-L A^{-1} B,
		\end{align*}
		is $L_\phi$-smooth and $m_\phi$-strongly convex, and $\phi_y$ is $L_{\phi_y}$-smooth. Finally, there exists $T_e\in[0,T)$ such that $\mu([0,T_e]) > 0$.
	\end{assumption}
	The reduced order SHDS associated with $\mathcal{H}_{\varepsilon}$ is the SHDS $\mathcal{H}_0$ defined by \eqref{SPSHDSreduced},
	where the flow map $\tilde{F}$ and the jump maps $\tilde{G}$ are defined by
	\begin{align*}
		\tilde{F}({x})&:=\{(p,-\beta p -\nabla \phi_u(u)-K\nabla \phi_y(Hu+d),1)\}, \\
		\tilde{G}({x},v)&:=\{(u,s p,v)~|~s\in[0,\rho]\}.
	\end{align*}
	Taking $K = H^\top$, we see that $\nabla \phi_u(u)+K\nabla \phi_y(d+Hu) = \nabla \phi(u)$,
	where $\phi(u):=\phi_u(u) + \phi_y(d+Hu)$, and that
	\begin{align*}
		\tilde{F}({x})&:=\{(p,-\beta p -\nabla \phi(u),1)\}.
	\end{align*}
	%
	%
\begin{proposition}\label{thm:heavy_ball_resets}
		Suppose that Assumption~\ref{asmp:stable_plant} holds, and let \tcb{$u^\star$
be the unique minimizer of $\phi$ over $\mathbb R^n$}. Define
$\mathcal A:=\{u^\star\}\times\{0\}\times[0,T]$ and
$\tilde{\mathcal A}:=\mathcal A\times\{-A^{-1}Bu^\star\}$.
Then there exists $\varepsilon^*>0$ such that, for every
$\varepsilon\in(0,\varepsilon^*)$ and with $K=H^\top$, the SP-SHDS
$\mathcal H_\varepsilon$ renders $\tilde{\mathcal A}$ UGASp.
	\end{proposition}

    The following example illustrates  Proposition 2.
\begin{example}\normalfont\label{example3}
		Let $u\in\mathbb{R}^2$, $z\in\mathbb{R}^2$, $T=100$, and $A=[-1,0;0,-1]$, $B=[1,0;0,1]$, $H=[1,0;0,1]$, $d=[1;-1]$.
		%
		%
		In addition, suppose that $\phi_u(u) = u^\top u,~\phi_y(y) = y^\top y$ and that $\mu$ is the (reflected) truncated Exponential distribution over $[0,T]$, i.e., 
		\begin{align*}
			\mu(dv) = \begin{cases}
				\frac{\mathrm{e}^{-(T-v)}}{1-\mathrm{e}^{-T}}\, dv& v\in[0,T] \\
				0 & \text{otherwise.}                
			\end{cases}
		\end{align*}
		The distribution $\mu$ models stochastic resetting of the momentum with (truncated) exponentially distributed waiting times. It can be shown that Assumption \ref{asmp:stable_plant} holds, and, therefore, for sufficiently small $\varepsilon$, system $\mathcal{H}_{\varepsilon}$ is UGASp. Figure~\ref{fig:exmp_2} shows a sample path for $\varepsilon=0.1$.  \QEDB 
		
		\begin{figure}
			    \centering
			    \includegraphics[width=0.95\linewidth]{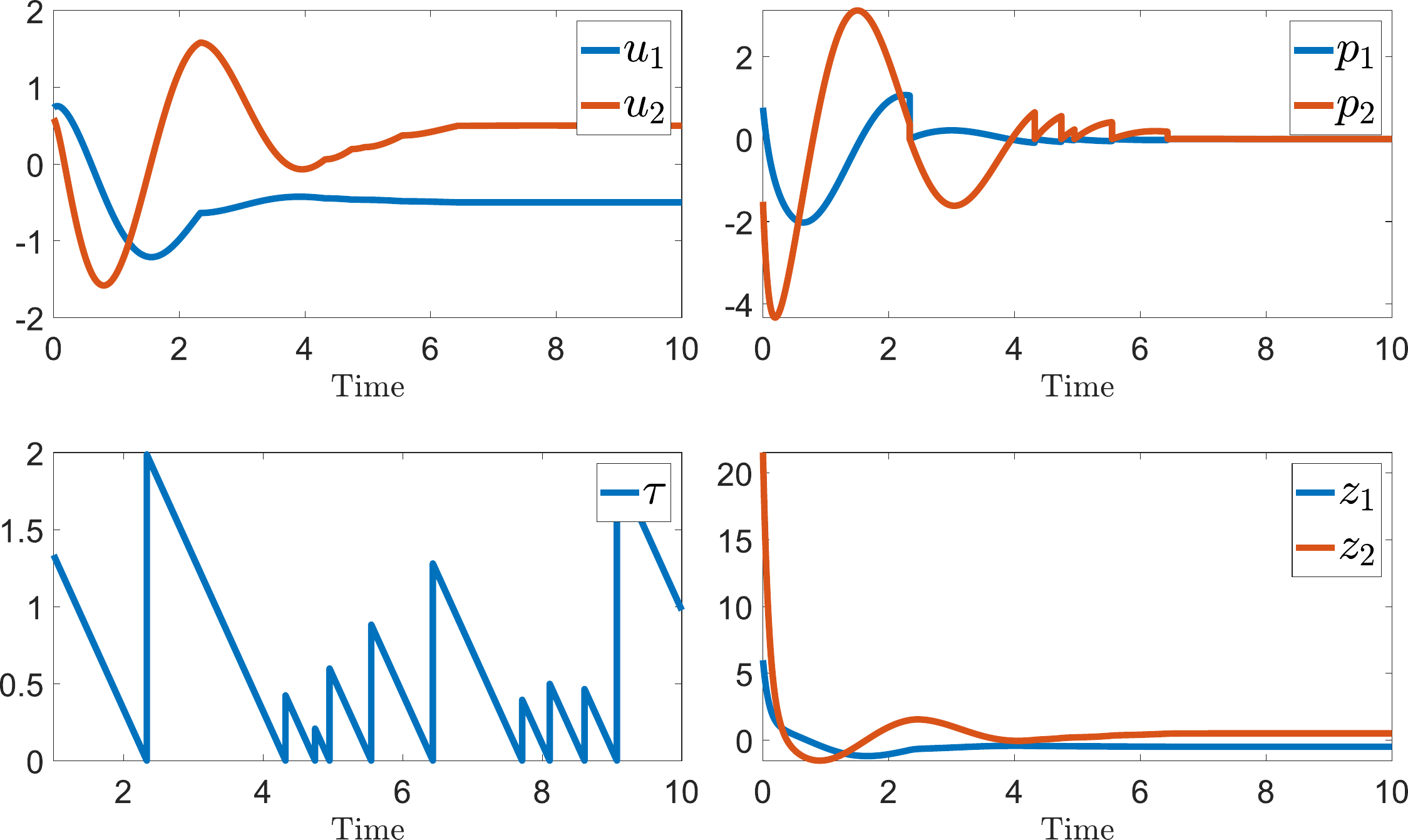}
			    \caption{Simulation results for feedback optimization problem of Example \ref{example3}.}
			    \label{fig:exmp_2}
			\end{figure}
	\end{example}
    
	\textbf{Proof of Proposition 2:}
	Define the function $W$ by
	\begin{align*}
		W(x,z):=(z+A^{-1}Bu)^\top P(z+A^{-1}Bu),
	\end{align*}
	where $P$ is the matrix from Assumption \ref{asmp:stable_plant},
	as well as the function $V$ by $V(x)=\frac{1}{2}p^\top p + \phi(u)-\phi(u^\star)$. Clearly, the function $W$ is radially unbounded and positive definite with respect to the slow manifold $\mathcal{M}(x):=\{z\in\mathbb{R}^m~|~z+A^{-1}B u = 0\}$.
	%
	%
	Indeed, since $|z|_{\mathcal{M}(x)} = \|z+A^{-1}B u\|$, we have that $\alpha_1(|z|_{\mathcal{M}(x)}) \leq W(x,z)\leq \alpha_2(|z|_{\mathcal{M}(x)}),$
	where the functions $\alpha_1$ and $\alpha_2$, defined by $\alpha_1(r):= \lambda_{\min}(P) r^2$ and $\alpha_2(r):= \lambda_{\max}(P) r^2$
	%
	%
	are both class $\mathcal{K}_\infty$ functions. 
	In addition, direct computation shows that 
	the function $W$ satisfies the inequality
	\begin{align*}
		\max_{h\in\partial_z W(x,z)} \langle h,f_z\rangle &\leq -k_z \varphi_z(|z|_{\mathcal{M}(x)})^2,
	\end{align*}
	for each $(x,z)\in C$ and for all $f_z\in F_z(x,z)$, where $k_z$ and $\varphi_z$ are defined by
	\begin{align*}
		k_z&:= - \lambda_{\max}(A^\top P + P A ), & \varphi_z(r)&:= |r|.
	\end{align*}
	On the other hand, since the function $\phi$ is strongly convex and $L_\phi$-smooth, it is clear that $V$ is positive definite and radially unbounded with respect to the compact set $\mathcal{A}:=\{u^\star\}\times\{0\}\times[0,T]$. Indeed, since $|x|_{\mathcal{A}} = \|(u-u^\star,p)\|$, we have that $\alpha_3(|x|_{\mathcal{A}})\leq V(x)\leq \alpha_4(|x|_{\mathcal{A}}),$
	%
	%
	for all $x\in C_x\cup D_x$, where $\alpha_3$ and $\alpha_4$, defined by $\alpha_3(r):=\frac{1}{2}\min\left\{1,m_\phi\right\} r^2$ and $\alpha_4(r):=\frac{1}{2}\max\left\{1,L_\phi\right\} r^2$
	%
	%
	are both $\mathcal{K}_\infty$ functions. In addition, for each $x\in C_x$ and all $f_x\in \tilde{F}(x)$, we compute that
	\begin{align*}
		\max_{h\in\partial_x V(x)}\langle h,f_x\rangle = -\beta p^\top p.
	\end{align*}
	Therefore, the function $V$ satisfies the inequality $\max_{h\in\partial_x V(x)}\langle h,f_x\rangle \leq -k_x\varphi_x(x)^2\leq 0$, for each $x\in C_x$ and all $f_x\in \tilde{F}(x)$, where $k_x:=\beta$ and $\varphi_x(x):=\|p\|$. Moreover, for every $x\in D_x$, we have
	\begin{align*}
		\int_{\mathbb R}\sup_{g\in\tilde G(x,v)}V(g)\,\mu(\mathrm dv)-V(x)\leq \frac{1}{2}(\rho^2-1)p^\top p.
	\end{align*}
	Hence, $V$ satisfies $\int_{\mathbb R}\sup_{g\in\tilde G(x,v)}V(g)\,\mu(\mathrm dv)
\leq V(x)-c_x\rho_x(x)$,
	%
	%
	for all $x\in D_x$, where $c_x:=\frac{1}{2}(1-\rho^2)$ and $\rho_x(x):=\|p\|^2$.
	Furthermore, it can be shown that, 
	for each $({x},{z})\in C$ and every $f_{{x}}\in F_{x}({x},{z})$, there exists $\tilde{f}_{{x}} \in \tilde{F}({x})$ such that
	\begin{gather*}
		\max_{h\in \partial_{x} W({x},{z})} \langle h,f_{x}\rangle \leq k_1 \varphi_{x}({x}) \varphi_{z}(|{z}|_{\mathcal{M}({x})})+ k_2\varphi_{z}(|{z}|_{\mathcal{M}({x})})^2,\\
		\max_{h\in \partial_{x} V({x})} \langle h,f_{x} - \tilde{f}_{x}\rangle \leq k_3 \varphi_{x}({x})\varphi_{z}(|{z}|_{\mathcal{M}({x})}),
	\end{gather*}
	where $k_2=0$ and we choose
$k_1>2\sigma_{\max}(B^\top A^{-\top}P)$ and
$k_3>\sigma_{\max}(H)\sigma_{\max}(L)L_{\phi_y}$.
	%
	%
	Finally, for all $({x},{z})\in D$, we have 
	%
	%
	%
	\begin{align*}
		\int_{\mathbb R}\sup_{g\in G(x,z,v)}W(g)\,\mu(\mathrm dv)
       = W(x,z)+k_5\rho_5(x),
	\end{align*}
	with $k_5:=0$ and any non-negative function $\rho_5({x})$, e.g., $\rho_5({x}):=\rho_{x}({x})$. By combining all of the above, it follows that the SP-SHDS $\mathcal{H}_{\varepsilon}$ satisfies Assumptions \ref{manifold},\ref{Assumption2},\ref{Assumption3},\ref{assumptionjumpsreduced},\ref{assumption6}, and \ref{jumpcoupledcondition} in Section \ref{section3}. Therefore, using $E_{\theta^*}(x,z)=(1-\theta^*) V (x) + \theta^* W(x,z)$ a computation as in the proof of Theorem \ref{theorem2} shows that there exists $\lambda > 0$ such that:

    \vspace{-0.9cm}
	\begin{align*}
		\max_{h\in \partial E_{\theta^*}(x,z)} \langle h,f\rangle \leq - \lambda (\varphi_x(x)^2 + \varphi_z(|z|_{\mathcal{M}(x)})^2)
	\end{align*}

    \vspace{-0.5cm}\noindent 
	for all $(x,z)\in C$, all $f\in F_x(x,z)\times\varepsilon^{-1}F_z(x,z)$, and all $\varepsilon\in (0,\varepsilon^*)$, and that

    \vspace{-0.8cm}
	\begin{align*}
		\int_{\mathbb R}\sup_{g\in G(x,z,v)} E_{\theta^*}(g)\,\mu(dv)
\leq E_{\theta^*}(x,z),
	\end{align*}

    \vspace{-0.5cm}\noindent 
	for all $(x,z)\in D$. A solution confined to a positive
level set of $E_{\theta^*}$ must satisfy $p=0$ and
$z=-A^{-1}Bu$ on every positive-duration flow interval.
On such an interval, $\dot p=-\nabla\phi(u)=0$, so
$u=u^\star$ and $E_{\theta^*}=0$, a contradiction. A complete solution confined to a positive level set
would therefore have to jump infinitely often without
flowing. Since $\mu([0,T_e])>0$ and $T_e<T$, almost surely an infinite sequence of resets includes one
with $v\leq T_e$, which requires flow before the next jump. This verifies item~(c) of Theorem~\ref{theorem2} and concludes the proof. \hfill$\blacksquare$ 
	%
	%
	%

    \vspace{-0.1cm}
	\subsection{Feedback Optimization with Spontaneously Switching Plant in the Loop}
	Let ${z}:=(\xi,q,\tau)\in\mathbb{R}^n\times\mathcal{Q}\times[0,T]$, where $\mathcal{Q}:=\{1,\dots,N\}$ and $T>0$, and let ${x}\in\mathbb{R}^{n_1}$. Let $v:=(v_1,v_2)\in\mathbb{R}^2$ and consider the SP-SHDS defined by \eqref{SPSHDS1}
	where the maps $F_{x}$, $F_{{z}}$, and $G$ are 
	\begin{subequations}
		\begin{align}
			F_{{x}}({x},{z})&:= \{-\nabla\phi_x(x)-H^\top\nabla\phi_y(L\xi+d)\},\\
			F_{{z}}({x},{z})&:=\{(A_q (\xi + B {x}), 0,s)~|~s\in[-\eta,0]\},\\
			G({x},z,v)&:=\{x\}\times\{(\xi,v_1,v_2)\},
		\end{align}
		for some $\eta\geq0$ and with $H = -LB$, and the flow set $C:=C_x\times C_z$ and the jump set $D:=D_x\times D_z$ are 
		\begin{align}
			C_x&:=\mathbb{R}^{n_1}, & C_z&:=\mathbb{R}^n\times\mathcal{Q}\times[0,T],\\
			D_x&:=\mathbb{R}^{n_1}, & D_z&:=\mathbb{R}^n\times\mathcal{Q}\times\{0\},
		\end{align}
	\end{subequations}
	for some $T>0$, and the definition of the measure $\mu$ coincides with that introduced in subsection \ref{subsec1}.
	The SP-SHDS $\mathcal{H}_{\varepsilon}$ models an interconnection between a gradient-based feedback optimization algorithm with a fast plant whose evolution is governed by a switching linear system with ``spontaneous" mode transitions. Associated with $\mathcal{H}_{\varepsilon}$ is the reduced order HDS defined by \eqref{SPSHDSreduced}
	where the flow map $\tilde{F}$ and the jump maps $\tilde{G}$ are
	\begin{align}
		\tilde{F}({x})&:=\{-\nabla\phi_x(x)- H^\top \nabla\phi_y(Hx+d)\}, \\
		\tilde{G}({x},v)&:=\{x\}.
	\end{align}
	We impose the following assumption to ensure that the boundary layer model as well as the reduced order SHDS are well-behaved (see Remark \ref{remark1}).
	\begin{assumption}\label{asmp:hurwitz_Aq}
		There exists $\sigma>0$, and a family of symmetric positive \tcb{definite} matrices $\{P_q\}_{q\in\mathcal{Q}}$, such that
		\begin{subequations}\label{eq:matrix_inequalities_2}
			\begin{align}
				A_q^\top P_q + P_q A_q + \sigma \eta T^{-1} P_q &\prec 0, \\
				\sigma^{-1}\log(1+\sigma)P- P_q &\prec 0,
			\end{align}
		\end{subequations}
		for all $q\in\mathcal{Q}$, where $P=\sum_{q\in\mathcal{Q}}\lambda_q P_q$. In addition, there exists  $m_{\phi},L_{\phi}$,$L_{\phi_y}>0$, such that the function
		\begin{align}
			\phi(x):=\phi_x(x) + \phi_y(H x + d),
		\end{align}
		is $L_\phi$-smooth and $m_\phi$-strongly convex, and $\phi_y$ is $L_{\phi_y}$-smooth.
	\end{assumption}
	Under Assumption \ref{asmp:hurwitz_Aq}, we have the following result.
	\begin{proposition}\label{proposition3}
		Suppose Assumption~\ref{asmp:hurwitz_Aq} holds, \tcb{and let $x^\star$
be the unique minimizer of $\phi$ over $\mathbb R^{n_1}$.}
Define $\mathcal A:=\{x^\star\}$ and
$\tilde{\mathcal A}:=\mathcal A\times\{-Bx^\star\}\times\mathcal Q\times[0,T]$.
Then there exists $\varepsilon^*>0$ such that, for every
$\varepsilon\in(0,\varepsilon^*)$, the SP-SHDS $\mathcal H_\varepsilon$
renders $\tilde{\mathcal A}$ UGASp.
	\end{proposition}
	The following example illustrates Proposition \ref{proposition3}.

    \begin{example}\normalfont\label{ex:switching_plant}
		Let $x,\xi\in\mathbb{R}^2$, $\mathcal{Q}=\{1,2\}$, $T=2$, and $A_1=[-1,3;0,-1]$, $A_2=[-1,0;-3,-1]$, $B=[-1,0;0,-1]$, $L=[1,0;0,1]$, $d=[1;-1]$, and  $\lambda_1=\lambda_2=1/2$.
		As in Example~2, the pair admits no common quadratic
        Lyapunov function with strict decrease. 
		For $\sigma=18$, direct computation shows that $P_1=[0.50, 0.75;0.75,2.75]$, $P_2=[2.75,-0.75;-0.75,0.50]$
		%
		satisfy the second inequality in \eqref{eq:matrix_inequalities_2}. Moreover, we have that $A_1^\top P_1 + P_1 A_1= A_2^\top P_2 + P_2 A_2 = [-1, 0; 0 , -1]\prec 0$. Consequently, for any $\eta<0.035$, the first inequality
in \eqref{eq:matrix_inequalities_2} is satisfied. In addition,
we take $\phi_x(x)=\tfrac14 x^\top x$ and
$\phi_y(y)=\tfrac14 y^\top y$.
        Figure \ref{fig:exmp_3} shows a  typical sample path obtained when $\varepsilon=0.01$. \hfill \QEDB

		\begin{figure}
			    \centering
			    \includegraphics[width=0.95\linewidth]{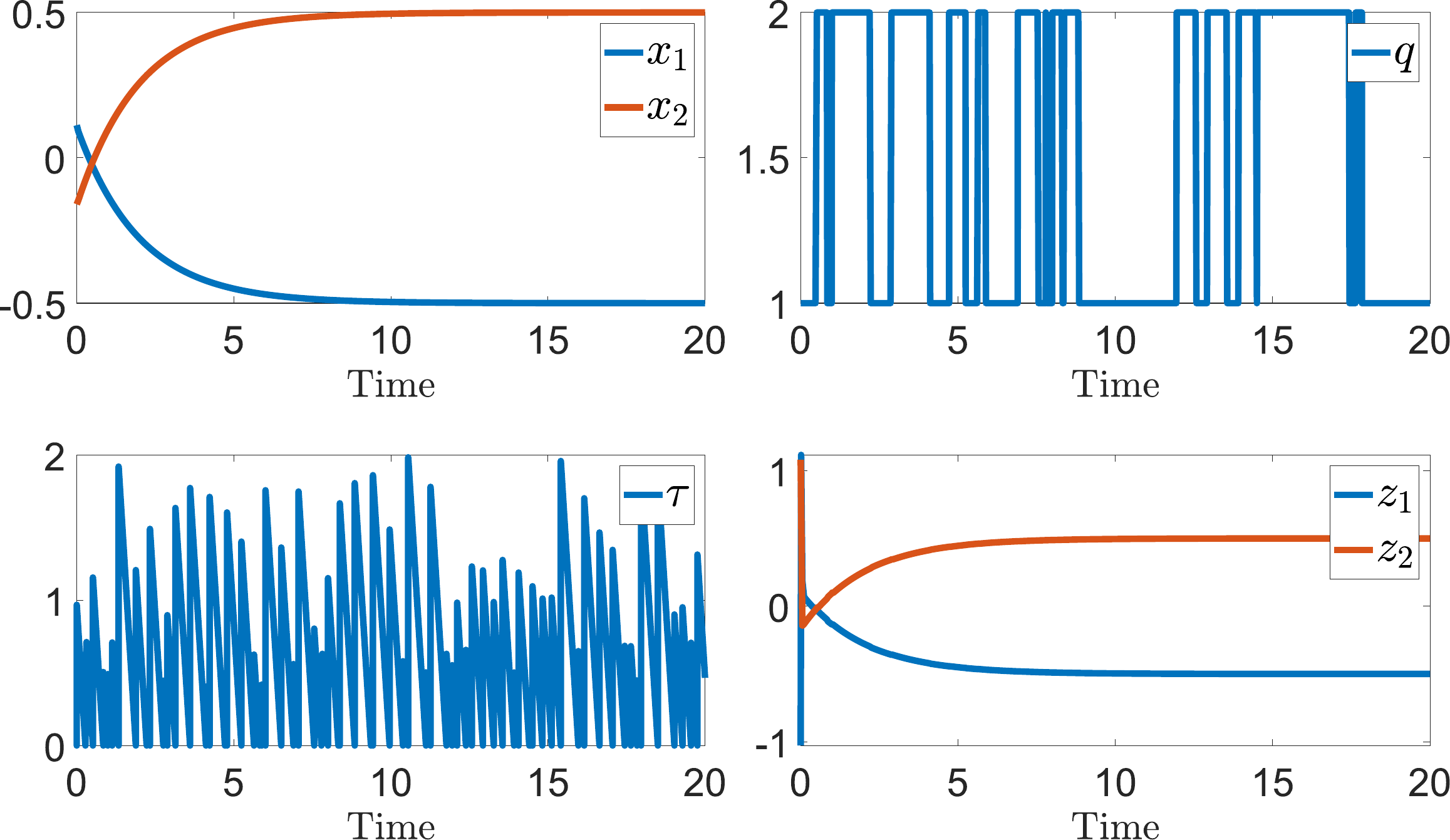}
			    \caption{Simulation results for Example~\ref{ex:switching_plant}.}
			    \label{fig:exmp_3}
			\end{figure}
	\end{example}
		
		\textbf{Proof of Proposition 3:} Define the function
		\begin{align*}
			W({x},{z}):=(\sigma\tau T^{-1}+1)^{-1} ({\xi}+Bx)^\top P_q({\xi}+Bx),
		\end{align*}
		and the function $V({x}):=\phi(x)-\phi(x^\star)$. Clearly, $W$ is positive definite and radially unbounded with respect to the slow manifold $\mathcal M(x)=\{-Bx\}\times\mathcal Q\times[0,T]$. Indeed, since $|{z}|_{\mathcal{M}({x})} = \|{\xi} +B x \|$, we have that $\alpha_1(|{z}|_{\mathcal{M}({x})})\leq W({x},{z})\leq \alpha_2(|{z}|_{\mathcal{M}({x})})$, where the functions $\alpha_1$ and $\alpha_2$, given by $\alpha_1(r):=(\sigma+1)^{-1}\min_{q\in\mathcal{Q}}\lambda_{\min}(P_{q}) r^2$, $\alpha_2(r):=\max_{q\in\mathcal{Q}}\lambda_{\max}(P_{q}) r^2$ are both $\mathcal{K}_\infty$ functions.  In addition, for each $({x},{z})\in C$ and for all $f_{{z}}\in F_{z}({x},{z})$,
		\begin{align}
			\max_{h\in \partial_{z} W({x},{z})}&\langle h,f_{{z}}\rangle \leq \frac{({\xi}+Bx)^\top Q_q ({\xi}+Bx)}{\sigma\tau T^{-1}+1}, \\
			Q_q&:=(A_q^\top P_{q} + P_{{q}}A_q + \sigma\eta T^{-1}P_q).
		\end{align}
		Hence, under Assumption \ref{asmp:hurwitz_Aq}, $W$ satisfies 
		\begin{align}
			\max_{h\in \partial_{z} W({x},{z})}\langle h,f_{{z}}\rangle &\leq - k_{z} \varphi_{z}(|{z}|_{\mathcal{M}({x})})^2,
		\end{align}
		for each $({x},{z})\in C$ and all $f_{{z}}\in F_{z}({x},{z})$, where $\varphi_{z}(r) = |r|$, and $k_{z}$ is defined by $k_{z}:= -(\sigma+1)^{-1}\max_{q\in\mathcal{Q}}\lambda_{\max}(A_q^\top P_{q} + P_{{q}}A_q + \sigma\eta T^{-1}P_q)$. Under Assumption \ref{asmp:hurwitz_Aq}, the constant $k_{z}$ is strictly positive. Moreover, direct computation gives that, for every $(x,z)\in D$, the function $W$ satisfies 
		\begin{align*}
			\int_{\mathbb{R}^2}&\sup_{g\in G({x},z,v)}W(g) ~\mu(\mathrm{d}v)-W({x},z) \\
			&=  ({\xi}+Bx)^\top(\sigma^{-1}\log(1+\sigma)P - P_q) ({\xi}+Bx).
		\end{align*}
		Therefore, the function $W$ satisfies the inequality
		\tcb{\begin{align*}			\int_{\mathbb{R}^2}\sup_{g\in G(x,z,v)} W(g)\,\mu(\mathrm{d}v)
        &\leq W(x,z)-c_z\rho_z\bigl(|z|_{\mathcal M(x)}\bigr),
		\end{align*}}
		for all $(x,z)\in D$, where $c_{z}:=- \max_{q\in\mathcal{Q}}\lambda_{\max}(\sigma^{-1}\log(1+\sigma)P- P_q)$ and $\rho_z(r):=r^2$.
		Note that, under Assumption \ref{asmp:hurwitz_Aq}, the constant $c_{z}$ is also strictly positive. 

		On the other hand, it is clear that $V$ is also positive definite and radially unbounded with respect to the compact set $\mathcal{A}$. Indeed, from Assumption \ref{asmp:hurwitz_Aq}, we have that $\alpha_3(|{x}|_{\mathcal{A}})\leq V({x}) \leq \alpha_4(|{x}|_{\mathcal{A}})$, for all ${x}\in C_x$, where the functions $\alpha_3(r):=m_\phi r^2/2$ and $\alpha_4(r):=L_\phi r^2/2$
		are both $\mathcal{K}_\infty$ functions.
		In addition, for each ${x}\in C_x$ and all $f_{x}\in \tilde{F}({x})$, we compute that
		\begin{align}
			\max_{h\in \partial_{x} V({x})}\langle h,f_{x}\rangle &\leq -\|\nabla\phi(x)\|^2 \leq -m_\phi^2|x|_{\mathcal A}^2,
		\end{align}
		where the second inequality follows from strong convexity. Therefore, under Assumption \ref{asmp:hurwitz_Aq}, the function $V$ satisfies the inequality $\max_{h\in \partial_{x} V({x})}\langle h,f_{x}\rangle\leq - k_{x} \varphi_{x} ({x})^2 \leq 0$ for each ${x}\in C_x$ and all $f_{x}\in \tilde{F}({x})$, where $k_{x}:=m_\phi^2$ and $\varphi_{x}({x}) = |{x}|_{\mathcal{A}}$.
		Note that, under Assumption \ref{asmp:hurwitz_Aq}, the constant $k_{x}$ is also strictly positive. Moreover, it can be shown that
		%
		for each $({x},{z})\in C$ and every $f_{{x}}\in F_{x}({x},{z})$, there exists $\tilde{f}_{{x}} \in \tilde{F}({x})$ such that
		\begin{gather*}
			\max_{h\in \partial_{x} W({x},{z})} \langle h,f_{x}\rangle \leq k_1 \varphi_{x}({x}) \varphi_{z}(|{z}|_{\mathcal{M}({x})})+ k_2\varphi_{z}(|{z}|_{\mathcal{M}({x})})^2,\\
			\max_{h\in \partial_{x} V({x})} \langle h,f_{x} - \tilde{f}_{x}\rangle \leq k_3 \varphi_{x}({x})\varphi_{z}(|{z}|_{\mathcal{M}({x})}),
		\end{gather*}
		where, with $P_{\max}:=\max_{q\in\mathcal Q}\|P_q\|$,
we choose $k_1>2L_\phi\|B\|P_{\max}$,
$k_2=2L_{\phi_y}\|B\|\|L\|\|H\|P_{\max}$, $k_3>L_\phi L_{\phi_y}\|H\|\|L\|$,
and $\|\cdot\|$ denotes the spectral norm for matrices.
Since $\tilde G(x,v)=\{x\}$, Assumption~\ref{interconnection2}
holds with $k_6=0$ and $\rho_6=\rho_z$. Together with $c_z>0$,
this verifies Theorem~\ref{theorem2}(b)(2). For sufficiently small $\varepsilon$, the flow estimate in
Theorem~\ref{theorem2} is strict outside $\tilde{\mathcal A}$.
A complete solution confined to a positive level set of
$E_{\theta^*}$ must therefore jump infinitely often without
flowing. Almost surely, each timer reset is positive, and
another jump requires the timer to flow to zero. This
verifies item~(c), and Theorem~\ref{theorem2} gives the result. \hfill$\blacksquare$


	\subsection{Recurrence in Linear Systems with Bounded Stochastic Inputs}
	Let ${x}:=(\xi,u,\tau)\in\mathbb{R}^n\times \mathcal{U}\times[0,T]$, where $\mathcal{U}\subset\mathbb{R}^n$ is a non-empty compact set, and let ${z}\in\mathbb{R}^{n_2}$. Consider the SP-SHDS $\mathcal{H}_{\varepsilon}$ in \eqref{SPSHDS1}, where the maps $F_{x}$, $F_{{z}}$, and $G$ are 
	\begin{subequations}
		\begin{align*}
			F_{{x}}({x},{z})&:=\{(A\xi+Bz+u, 0,1)\},\\
			F_{{z}}({x},{z})&:= \{H\xi+Lz\}\\
			G({x},{z},v)&:=\{(\tcb{\xi},v,0)\}\times\{z\},
		\end{align*}
		and $C:=C_x\times C_z$ and $D:=D_x\times D_z$ are defined by
		\begin{align*}
			C_x&:=\mathbb{R}^n\times \mathcal{U}\times[0,T], & C_z&:= \mathbb{R}^{n_2}, \\
			D_x&:=\mathbb{R}^n\times \mathcal{U}\times\{T\}, & D_z&:=\mathbb{R}^{n_2},
		\end{align*}
	\end{subequations}
	for some $T>0$.
	In this model, $v\in\mathbb{R}^n$ is a random variable distributed according to some probability distribution $\mu$ with support on $\mathcal{U}$. 

    \vspace{-0.1cm}
 \begin{assumption}\label{asmp:stable_plant_2}
		\tcb{There exist $P\succ0$ and $Q\succ0$}, such that the matrices $L$ and $\tilde{A}=A-BL^{-1}H$ satisfy the inequalities $\tilde{A}^\top P + P \tilde{A} \prec 0$ and $L^\top Q + Q L \prec 0$.
	\end{assumption}
	The reduced order SHDS associated with $\mathcal{H}_{\varepsilon}$ is the SHDS $\mathcal{H}_0$ defined by \eqref{SPSHDSreduced}, where $\tilde{F}$ and $\tilde{G}$ are defined by
	\begin{align}
		\tilde{F}({x})&:=\{(\tilde{A}\xi + u,0,1)\}, &
		\tilde{G}({x},v)&:=\{(\xi,v,0)\}.
	\end{align}
For any fixed $\delta>0$, define
\begin{equation}\label{eq:recurr_set_44}
\mathcal O_x:=r\mathbb B^\circ
\times(\mathcal U+\delta\mathbb B^\circ)
\times(-\delta,T+\delta).
\end{equation}
	with $r:=\inf \{\tilde{r}\in\mathbb{R}_{>0} ~|~ \kappa (\xi) \leq0 , \, \forall \xi\in\mathbb{R}^n,\,\|\xi\|\geq \tilde{r}\}$, where $\kappa(\xi):=(1-\chi_1)\lambda_{\max}(\tilde{A}^\top P + P \tilde{A})\|\xi\|^2 + 2\sup_{u\in \mathcal{U}}\|P u\|\|\xi\| - \chi_1 \lambda_{\max}(\tilde{A}^\top P + P \tilde{A}) \tilde{c}_x^2$, for any $\chi_1\in(0,1)$ and any $\tilde{c}_x>0$.
	Under Assumption \ref{asmp:stable_plant_2}, $r>0$ and the set $\mathcal{O}_x$ is bounded and open. 

    \vspace{-0.1cm}
    \begin{proposition}\label{prop:bounded_input_recurrence}.
	Suppose Assumption~\ref{asmp:stable_plant_2} holds.
Then, for any $\chi,\tilde c_x>0$ and $\chi_1\in(0,1)$, there
exists $\varepsilon^*>0$ such that, for every
$\varepsilon\in(0,\varepsilon^*)$, the set $\mathcal O_\chi$
in \eqref{eq:recurr_set_331}, with $\mathcal O_x$ given by
\eqref{eq:recurr_set_44}, is UGR for $\mathcal H_\varepsilon$
relative to $\varpi(y)=\|y\|$. 
	\end{proposition}
\vspace{-0.2cm}
\begin{example}\normalfont\label{ex:bounded_input_recurrence}
		Let $\xi\in\mathbb{R}^2$, $z\in\mathbb{R}^2$, $\mathcal{U}=\mathbb{B}\subset\mathbb{R}^2$, $T=1$, and $A=[-2,2;-1,0]$, $B=[0,1;1,0]$, $L=[-1,0;0,-1]$, $H=[1,-1;1,1]$.
		Direct computation gives $\tilde{A}=[-1,3;0,-1]$,
		%
		%
		which is Hurwitz. We consider the uniform distribution over $\mathcal{U}$, and we show in Figure~\ref{fig:exmp_4} a sample path for $\varepsilon=0.1$. \hfill \QEDB

		\begin{figure}[t!]
			    \centering
			    \includegraphics[width=0.95\linewidth]{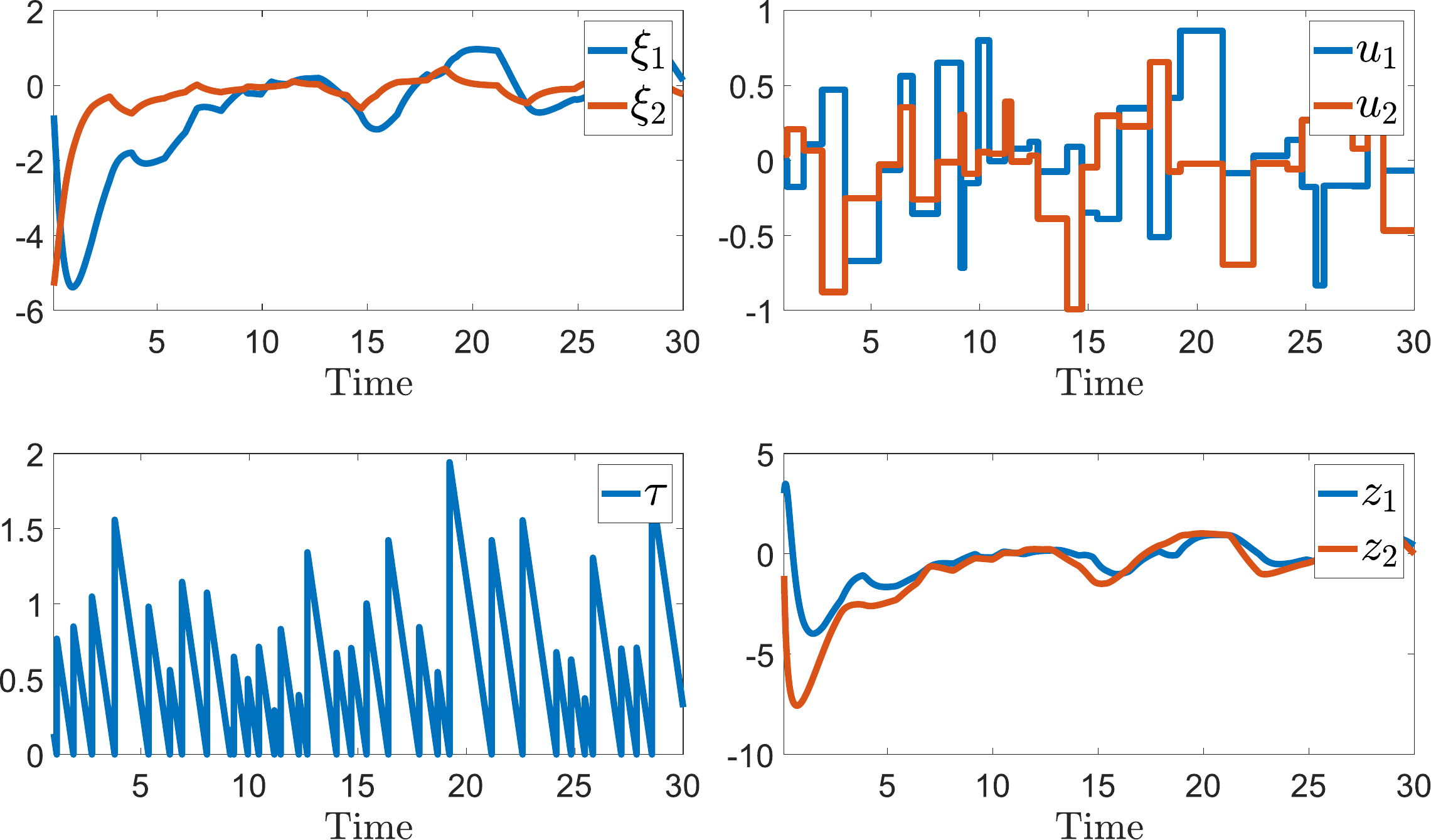}
			    \caption{Simulation results for Example~\ref{ex:bounded_input_recurrence}.}
			    \label{fig:exmp_4}
			\end{figure}
	\end{example}

    \vspace{-0.2cm}
	\textbf{Proof:} Let $V$ and $W$ be given by
	\begin{align*}
		V(x)&:=\xi^\top P \xi , \\
		 W(x,z)&:=(z+L^{-1}H\xi)^\top Q (z+L^{-1}H\xi),
	\end{align*}
	where $P$ and $Q$ are as given in Assumption \ref{asmp:stable_plant_2}.
	Clearly, the function $W$ is radially unbounded and positive definite with respect to the slow manifold $\mathcal{M}(x):=\{z\in\mathbb{R}^{n_2}~|~ z+L^{-1}H\xi= 0\}$. Indeed, since $|z|_{\mathcal{M}(x)} = \|z+L^{-1}H\xi\|$, we have that $\alpha_1(|z|_{\mathcal{M}(x)}) \leq W(x,z)\leq \alpha_2(|z|_{\mathcal{M}(x)})$, where the functions $\alpha_1$ and $\alpha_2$, defined by $\alpha_1(r):= \lambda_{\min}(Q) r^2$
		$\alpha_2(r)= \lambda_{\max}(Q) r^2$, are both class $\mathcal{K}_\infty$ functions. Also, for each $(x,z)\in C$ and $f_z\in F_z(x,z)$ we have
	\begin{align}
		\max_{h\in\partial_z W(x,z)} \langle h,f_z\rangle &\leq -k_z \varphi_z(|z|_{\mathcal{M}(x)})^2,
	\end{align}
	where $k_z$ and $\varphi_z$ are defined by $k_z:= - \lambda_{\max}(L^\top Q + Q L)$ and $\varphi_z(r):= |r|$. On the other hand, for all $f_x\in \tilde{F}(x)$ with $x\in C_x$, and by adding and subtracting terms, we compute that
	\begin{align*}
		\max_{h\in\partial_x V(x)}&\langle h,f_x\rangle = \xi^\top(\tilde{A}^\top P + P \tilde{A})\xi + 2 \xi^\top P u \\
		&\leq \lambda_{\max}(\tilde{A}^\top P + P \tilde{A})\|\xi\|^2+ 2\sup_{u\in \mathcal{U}}\|P u\|\|\xi\| \\
		&\leq \chi_1 \lambda_{\max}(\tilde{A}^\top P + P \tilde{A})(\|\xi\|^2 + \tilde{c}_x^2) + \kappa(\xi).
	\end{align*}
	Thus, there exists $\nu>0$ such that 
	\begin{align}
		\max_{h\in\partial_x V(x)}\langle h,f_x\rangle \leq -k_x\varphi_x(x)^2 + \nu\mathbb{I}_{\mathcal{O}_x}(x) ,
	\end{align}
	for all $f_x\in \tilde{F}(x)$ with $x\in C_x$, where $k_x:=-\chi_1\lambda_{\max}(\tilde{A}^\top P + P \tilde{A})$, $\varphi_x(x):= \sqrt{\|\xi\|^2+\tilde{c}_x^2}$. Furthermore, for all $f_{{x}}\in F_{x}({x},{z})$ with $({x},{z})\in C$, we get
	\begin{align*}
		&\max_{h\in \partial_{x} W({x},{z})} \langle h,f_{x}\rangle \leq 2 \sigma_{\max}(H^\top L^{-\top}Q) |{z}|_{\mathcal{M}({x})}\|\tilde{A}\xi+u\|\\
		&+2 \lambda_{\max}(\mathrm{Sym}(B^\top H^\top L^{-\top}Q)) |{z}|_{\mathcal{M}({x})}^2\\
		&\leq 2 \sigma_{\max}(H^\top L^{-\top}Q) |{z}|_{\mathcal{M}({x})}(\sigma_{\max}(\tilde{A})\|\xi\|+\sup_{u\in\mathcal{U}}\|u\|) \\
		&+2 \lambda_{\max}(\mathrm{Sym}(B^\top H^\top L^{-\top}Q)) |{z}|_{\mathcal{M}({x})}^2,
	\end{align*}
	and that, with $\tilde{f}_{{x}} = f_{{x}}+(\tilde{A} \xi-(A \xi + Bz),0,0) \in \tilde{F}({x})$, 
	\begin{align*}
		\max_{h\in \partial_{x} V({x})} \langle h,f_{x} - \tilde{f}_{x}&\rangle =\langle 2 P \xi ,B ({z}+L^{-1}H \xi)\rangle \\
		&\leq 2\sigma_{\max}(P) \sigma_{\max}(B) \|\xi\| |{z}|_{\mathcal{M}({x})}.
	\end{align*}
    Using $\|\xi\|\leq\varphi_x(x)$, choose $k_1,k_2,k_4$
to bound the coefficients of $\|\xi\|s$, $s^2$, and $s$,
respectively, in the bound for $W$, and choose $k_3>0$
to bound the coefficient of $\|\xi\|s$ in the bound for $V$,
where $s:=|z|_{\mathcal M(x)}$.
	%
    Set $\mathcal A=\overline{\mathcal O_x}$ and
$\varpi_x(x)=\|\xi\|$. Since $(0,u,\tau)\in\mathcal A$
on admissible slow states, $|x|_{\mathcal A}\leq\|\xi\|$,
which gives the comparison bounds on $V$.
The map $\mathcal M(x)=\{-L^{-1}H\xi\}$ satisfies
Assumption~\ref{manifold}, and the preceding estimates
verify item~(a) of Theorem~\ref{theorem4}, with
$m_\chi=\chi>0$. Jumps preserve $\xi$ and $z$, hence
$E_{\theta^*}$, giving item~(b). For sufficiently small
$\varepsilon$, $E_{\theta^*}$ strictly decreases along
flows outside $\mathcal O_\chi$. Since successive jumps
are separated by $T>0$, every complete solution has a
positive-duration flow interval. This verifies item~(c),
and Theorem~\ref{theorem4} yields the result.
\hfill$\blacksquare$

     \vspace{-0.2cm}
	\section{Conclusions}
	\label{sec_conclusions}

    \vspace{-0.3cm}
This paper studies Foster-based sufficient conditions for uniform global asymptotic stability in probability and uniform global recurrence in a class of stochastic hybrid dynamical systems with multi-time scale continuous dynamics. The framework relies on composite Lagrange–Foster and Lyapunov–Foster functions constructed from reduced and boundary-layer subsystems. As illustrated via different examples and applications, the resulting tools support the design of multi-time-scale feedback controllers that explicitly account for hybrid and stochastic effects. \tcb{Future research directions will focus on incorporating multi-time scale stochastic flows.}

        \vspace{-0.3cm}
\bibliographystyle{unsrt} 
{\small
	\bibliography{dissertationA.bib}
}

     \clearpage
 	\section*{Appendix}

 	The following auxiliary Lemmas correspond to \cite[Corollary 2 and 3, pp. 39-40]{Clarke}.
 	\begin{lemma}\label{lemma1}
 		Let $f_1,f_2$ be locally Lipschitz and let $s_1,s_2\in\mathbb{R}$. Then $\partial (s_1f_1+s_2f_2)(x)\subset s_1\partial f_1(x)+s_2\partial f_2(x)$, for all $x\in\text{dom} f_1\cap \text{dom} f_2$. Moreover, if $f_1,f_2$ are regular and $s_1,s_2\geq0$, then $\partial (s_1f_1+s_2f_2)(x)=s_1\partial f_1(x)+s_2\partial f_2(x)$.
 	\end{lemma}
 	The next Lemma is given by  \cite[Corollary 3, pp. 40]{Clarke}.
 	\begin{lemma}
 		Let $f(x,y)$ be locally Lipschitz and regular. Then $\partial f(x,y)\subset \partial_x f(x,y) \times \partial_y f(x,y)$, for all $(x,y)\in\text{dom}f$.
 	\end{lemma}
 	The next lemma gives some useful properties of the generalized gradient.
 	\begin{lemma}\cite[Prop. 1.5]{Clarke2}
 		Let $f:\mathbb{R}^n\to\mathbb{R}$ be locally Lipschitz. The Clarke generalized gradient of $f$, denoted $\partial f:\text{dom}~f\rightrightarrows\mathbb{R}^n$, is an outer-semicontinuous set-valued mapping with non-empty, compact and convex values.
 	\end{lemma}
 	\begin{lemma}
 		Let $f:\mathbb{R}^n\to\mathbb{R}$ be locally Lipschitz, and let $\phi:\text{dom}(\phi)\to\text{dom}(f)$ be an absolutely continuous function. Then, for almost all $t\in\text{dom}(\phi)$, there exists a $u\in\partial f(\phi(t))$ such that $\frac{d}{dt}(f\circ\phi)(t)=\langle u, \dot{\phi}(t)\rangle$.
 	\end{lemma}
 	Let $F:\mathbb{R}^n\rightrightarrows\mathbb{R}^n$, and let $f$ be a locally Lipschitz function. The set-valued Clarke derivative is defined as
 	\begin{equation}
 		\mathcal{C}_f(x;F):=\{\langle u,l\rangle :u\in \partial f(x),l\in F(x)\},
 	\end{equation}
 	The next Lemma is instrumental to analyze the behavior of trajectories of set-valued dynamical systems.
 	\begin{lemma}
 		Let $F:\mathbb{R}^n\rightrightarrows\mathbb{R}^n$. Given a locally Lipschitz function $V:\text{dom}(V)\to\mathbb{R}$, any solution $\phi$ to $\dot{x}\in F(x)$ satisfies $\tcb{\frac{d}{dt}(V\circ\phi)(t)\in \mathcal{C}_V(\phi(t);F)}$ for almost all $t\in\text{dom}(\phi)$.
 	\end{lemma}
 	The following lemmas make use of some of the previous results.
     \begin{lemma}
 Let $V:\mathbb R^n\to\mathbb R$ and
 $W:\mathbb R^n\times\mathbb R^m\to\mathbb R$ be locally
 Lipschitz, and define $E(x,z)=V(x)+W(x,z)$. Then, for every
 $f=(f_x,f_z)\in\mathbb R^n\times\mathbb R^m$,
 \[
 \max_{e\in\partial E(x,z)}\langle e,f\rangle
 \leq \max_{v\in\partial V(x)}\langle v,f_x\rangle
    + \max_{w\in\partial W(x,z)}\langle w,f\rangle.
 \]
 \end{lemma}
 \textbf{Proof:}
 The Clarke sum rule gives
 $\partial E(x,z)\subseteq
 (\partial V(x)\times\{0_m\})+\partial W(x,z)$.
 Taking the maximum of the inner product with $f$
 gives the result. \hfill$\blacksquare$

 	\vspace{0.1cm}
 	\begin{lemma}
 		Let $W:\mathbb{R}^n\times\mathbb{R}^m\to\mathbb{R}$ be locally Lipschitz and regular. Then, for any $(x,y)\in\mathbb{R}^n\times\mathbb{R}^m$, $f=(f_x,f_y)\in\mathbb{R}^n\times\mathbb{R}^m$, we have: 
 		\begin{equation*}
 			\max_{\substack{w\in \partial W(x,y)}}\langle w,f\rangle\leq \max_{\substack{w_x\in \partial_x W(x,y)}}\langle w_x,f_x\rangle+\max_{\substack{w_y\in \partial_y W(x,y)}}\langle w_y,f_y\rangle.
 		\end{equation*}
 	\end{lemma}
	
 	\vspace{-0.4cm}
 	\textbf{Proof:}  Using the definition of inner product
 	\begin{equation*}
 		\max_{\substack{w\in \partial W(x,y)}}\langle w,f\rangle=\max_{\substack{(w_x,w_y)\in \partial W(x,y)}}\Big(\langle w_x, f_x\rangle+\langle w_y,f_y\rangle\Big),
 	\end{equation*}
 	and since $\partial W(x,y)\subset \partial_x W(x,y)\times \partial_y W(x,y)$, the inequality follows directly.\hfill $\blacksquare$
 	\begin{lemma}
 		Let $f_x,\tilde{f}_x\in\mathbb{R}^n$, and $V:\mathbb{R}^n\to\mathbb{R}_{\geq0}$ be a regular Lipschitz function. Then, the following holds:
 		\begin{equation*}
 			\max_{\substack{v\in \partial V(x)}}\langle v,f_x \rangle-\max_{\substack{v\in \partial V(x)}}\langle v,\tilde{f}_x \rangle\leq \max_{\substack{v\in \partial V(x)}}\langle v,f_x-\tilde{f}_x \rangle 
 		\end{equation*}
 		for all $x\in\mathbb{R}^n$.
 	\end{lemma}
 	\textbf{Proof:} For each $x$, $\partial V(x)$ is convex and compact. We define the function $h_x:\mathbb{R}^n\to\mathbb{R}$ given by:
 	\begin{equation}
 		h_x(a):=\max_{v\in\partial V(x)}\langle v,a \rangle,
 	\end{equation}
 	which is the support function with respect to the set $\partial V(x)$. Since the support function of a compact and convex set is subadditive, it satisfies $h_x(a+c)\leq h_x(a)+h_x(c)$, which is equivalent to $h_x(a+c)-h_x(a)\leq h_x(c)$. Setting $c=f_x-\tilde{f}_x$ and $a=\tilde{f}_x$ leads to $$h_x(f_x)-h_x(\tilde{f}_x)\leq h_x(f_x-\tilde{f}_x),$$ which is the desired inequality. \hfill $\blacksquare$
 	\begin{lemma}\label{lemmaKbounds}
 \tcb{Suppose that $\tilde{\mathcal A}$ is nonempty, that
 Assumptions~\ref{manifold}, \ref{Assumption2}, and
 \ref{Assumption3} hold, and that
 $\alpha_i\in\mathcal K_\infty$ for $i=1,2,3,4$ and
 $\varpi_x(x)=|x|_{\mathcal A}$}. Then, for every fixed \(\theta\in(0,1)\), there exist \(\underline\alpha_\theta,\overline\alpha\in\mathcal K_\infty\) such that
 \begin{equation}\label{upperbound122}
 \underline\alpha_\theta\bigl(|y|_{\tilde{\mathcal A}}\bigr)
 \leq E_\theta(y)
 \leq\overline\alpha\bigl(|y|_{\tilde{\mathcal A}}\bigr)
 \end{equation}
 for all \(y\in \mathcal{S}:=C\cup D\cup G(D\times\mathcal V)\). \tcb{The upper comparison function can be chosen independently of \(\theta\).}     
 	\end{lemma}
	
 	\textbf{Proof:}  \tcb{For $y=(x,z)\in\mathcal S$, define $a:=|x|_{\mathcal A}$, $b:=|z|_{\mathcal M(x)}$, $d:=|y|_{\tilde{\mathcal A}}$. The assumed bounds on $V$ and $W$ imply the inequalities
 \begin{equation}\label{eq:lemma9_basic_bounds}
 \begin{aligned}
     (1-\theta)\alpha_3(a)+\theta\alpha_1(b)\leq E_{\theta}(y)\leq \alpha_4(a)+\alpha_2(b),
 \end{aligned}
 \end{equation}
 where $\alpha_i\in\mathcal K_{\infty}$ for
 $i\in\{1,2,3,4\}$. Under the given assumptions, note that the set $\tilde{\mathcal A}
     :=
     \{(r_1,r_2):r_1\in\mathcal A,\,
                      r_2\in\mathcal M(r_1)\}$ is nonempty and compact. 
 }

 \tcb{To establish the upper bound in \eqref{upperbound122}, choose $(x^*,z^*)\in\tilde{\mathcal A}$ such that $d^2=\|x-x^*\|^2+\|z-z^*\|^2$. Such a point exists by compactness of $\tilde{\mathcal A}$. Clearly, $a\leq\|x-x^*\|\leq d$. Moreover, the triangle inequality for distance and
 Assumption~\ref{manifold} give $b\leq\|z-z^*\|+\operatorname{dist}(z^*,\mathcal M(x))\leq
     \|z-z^*\|+\alpha(\|x-x^*\|)\leq d+\alpha(d)$. Consequently, \eqref{eq:lemma9_basic_bounds} implies
 \[
     E_{\theta}(y)
     \leq \alpha_4(d)+\alpha_2(d+\alpha(d))
     =:\overline{\alpha}(d).
 \]
 Since $\alpha\in\mathcal K$, the function
 $r\mapsto r+\alpha(r)$ belongs to $\mathcal K_{\infty}$.
 Thus, $\overline{\alpha}\in\mathcal K_{\infty}$.
 This function does not depend on $\theta$.}

 \tcb{To establish the lower bound in \eqref{upperbound122}, we first construct a comparison function relating
 $d$ to $a+b$.
 For $r\geq0$, define
 \[
     \omega(r):=
     \sup\left\{
        |(\xi,\zeta)|_{\tilde{\mathcal A}}:
        |\xi|_{\mathcal A}\leq r,\,
        \zeta\in\mathcal M(\xi)
     \right\}.
 \]
 The function $\omega$ is nondecreasing and finite
 for every finite $r$.
 Indeed, $\{\xi:|\xi|_{\mathcal A}\leq r\}$ is compact,
 and local boundedness of $\mathcal M$ implies
 boundedness of the corresponding pairs $(\xi,\zeta)$.
 Also, $\omega(0)=0$.}

 \tcb{We claim that
 \begin{equation}\label{eq:lemma9_omega_zero}   \lim_{r\downarrow0}\omega(r)=0.
 \end{equation}
 Otherwise, there would exist $\delta>0$ and sequences
 $r_k\downarrow0$, $\xi_k$, and $\zeta_k$ such that $|\xi_k|_{\mathcal A}\leq r_k$,
     $\zeta_k\in\mathcal M(\xi_k)$, $|(\xi_k,\zeta_k)|_{\tilde{\mathcal A}}
     \geq\delta$. Compactness of $\mathcal A$ and local boundedness
 of $\mathcal M$ allow extraction of a subsequence
 satisfying $\xi_k\to\xi^*\in\mathcal A,\qquad
     \zeta_k\to\zeta^*$. Outer semicontinuity of $\mathcal M$ yields
 $\zeta^*\in\mathcal M(\xi^*)$.
 Hence $(\xi^*,\zeta^*)\in\tilde{\mathcal A}$,
 contradicting the displayed distance bound.
 This proves \eqref{eq:lemma9_omega_zero}.}

 \tcb{These properties imply that $\omega$ admits a
 class-$\mathcal K_{\infty}$ majorant.
 For completeness, one can construct such a majorant
 as follows. For each integer $j$, set $b_j:=2^j+\omega(2^{j+1})$. Then $b_{j+1}>b_j$, $b_j\to0$ as $j\to-\infty$,
 and $b_j\to\infty$ as $j\to\infty$.
 Define $\beta(0):=0$ and let $\beta$ be the
 piecewise-affine function satisfying
 $\beta(2^j)=b_j$ for every integer $j$.
 It follows that $\beta\in\mathcal K_{\infty}$.
 Furthermore, if $r\in[2^j,2^{j+1}]$, then $\beta(r)\geq b_j
     \geq\omega(2^{j+1})
     \geq\omega(r)$. Thus $\omega(r)\leq\beta(r)$ for all $r\geq0$.}

 \tcb{Choose $z_0\in\mathcal M(x)$ such that
 $\|z-z_0\|=b$.
 The minimum is attained because $\mathcal M(x)$
 is nonempty and closed.
 Using the triangle inequality,
 we obtain
 \begin{align*}
     d&\leq\|z-z_0\|+
     |(x,z_0)|_{\tilde{\mathcal A}}\leq b+\omega(a)\leq b+\beta(a)\\
     &\leq (a+b)+\beta(a+b).
 \end{align*}
 For each fixed $\theta\in(0,1)$, define $q_{\theta}(r):=
     \min\{(1-\theta)\alpha_3(r),\,
                          \theta\alpha_1(r)\}$. Since both functions in the minimum belong to
 $\mathcal K_{\infty}$, so does $q_{\theta}$. Let $h(r):=r+\beta(r)$ and note that $h\in\mathcal K_{\infty}$, and therefore $h^{-1}(d)\leq a+b$. Then, using \eqref{eq:lemma9_basic_bounds}, we obtain
 \begin{align*}
     E_{\theta}(y)\geq q_{\theta}(a)+q_{\theta}(b)\geq q_{\theta}(\max\{a,b\})\geq q_{\theta}\left(\frac{a+b}{2}\right),
 \end{align*}
 which implies that $E_{\theta}(y)\geq q_{\theta}\left(\frac{h^{-1}(d)}{2}\right)$. Consequently, the function
 \[
     \underline{\alpha}_{\theta}(r)
     :=
     q_{\theta}\left(\frac{h^{-1}(r)}{2}\right)\in \mathcal K_{\infty}
 \]
 establishes the lower bound in \eqref{upperbound122}.}     \hfill $\blacksquare$

 \begin{lemma}\label{lemma10}
 \tcb{Under the hypotheses of Theorem 1, fix any \(\lambda>0\) and let \(p\) be defined as in \eqref{flowproof1}. There exists a continuous positive definite function
 $\rho_F:\mathbb R_{\geq0}\to\mathbb R_{\geq0}$ such that
 $\rho_F(|y|_{\tilde{\mathcal A}})\leq\lambda p(y)$ for all
 $y\in C$.}
 \end{lemma}

 \tcb{\textbf{Proof of Lemma \ref{lemma10}:}  For any $0<r\leq R<\infty$, the function $p$ is
 bounded away from zero on every nonempty set $\{y\in C:r\leq |y|_{\tilde{\mathcal A}}\leq R\}$. Indeed, if $p(y_k)\to0$ along a sequence in this set, then
 boundedness of $x_k$ and continuity and positive
 definiteness of $\varphi_x$ imply
 $|x_k|_{\mathcal A}\to0$. Also,
 $|z_k|_{\mathcal M(x_k)}$ is bounded by local boundedness
 of $\mathcal M$ and nonemptiness of $\mathcal M(x_k)$ on
 $C_x$. Since $\varphi_z$ is continuous and positive
 definite, it follows that $|z_k|_{\mathcal M(x_k)}\to0$.
 The distance comparison established in the lower-bound
 proof of Lemma~\ref{lemmaKbounds} then implies
 $|y_k|_{\tilde{\mathcal A}}\to0$, contradicting the lower
 bound $r$. For each $j\in\mathbb Z$, define $$a_j:=
     \min\left\{
         2^j,\,
         \inf_{\substack{y\in C\\
         2^{j-1}\leq |y|_{\tilde{\mathcal A}}\leq 2^{j+1}}}
         \lambda p(y)
     \right\},$$ where the infimum of an empty set is $+\infty$.
 The preceding argument gives $0<a_j\leq 2^j$.
 Set $\rho_F(0)=0$, set $\rho_F(2^j)=a_j$ for every
 $j\in\mathbb Z$, and interpolate linearly on each interval
 $[2^j,2^{j+1}]$. The resulting function is continuous
 and positive definite. If $|y|_{\tilde{\mathcal A}}$ lies
 in $[2^j,2^{j+1}]$, then both $a_j$ and $a_{j+1}$ are
 bounded above by $\lambda p(y)$. Thus,
 $\rho_F(|y|_{\tilde{\mathcal A}})\leq\lambda p(y)$.} \hfill $\blacksquare$

\end{document}